\documentclass[12pt]{amsart}
\usepackage{amscd,amssymb,verbatim}
\usepackage{amssymb,amsmath,amsthm,epsfig}

\newcommand{\N}{\mathbb{N}}

\newcommand{\R}{\mathbb{R}}

\newcommand{\supp}{{\rm supp}\,}

\theoremstyle{plain}
\newtheorem{Thm}{Theorem}[section]
\newtheorem{Prop}[Thm]{Proposition}
\newtheorem{Lem}[Thm]{Lemma}
\newtheorem{Cor}[Thm]{Corollary}
\newtheorem{Rmk}[Thm]{Remark}

\newtheorem{Qst}[Thm]{Question}

\newtheorem{Def}[Thm]{Definition}

\theoremstyle{remark}

\begin{document}

\title{On the ``Multiple of the Inclusion Plus Compact'' Problem}
\author{G. Androulakis, F. Sanacory$^*$}

\address{Department of Mathematics, University of South Carolina, Columbia, SC
29208.}
\email{giorgis@math.sc.edu, sanacory@math.sc.edu}
\date{January 12, 2007}

\thanks{$^*$The present paper is part of the Ph.D thesis of the second
author which will be prepared at the University of South Carolina under the
supervision of the first author}

\subjclass{46A32, 47B07}
\maketitle

\markboth{G. Androulakis, F. Sanacory}{On the ``Multiple of the Inclusion Plus Compact'' Problem}

\noindent
\begin{abstract} The ``multiple of the inclusion plus compact problem'' which
was posed by T.W.~Gowers in 1996 and Th.~Schlumprecht in
2003, asks whether for every infinite dimensional Banach space
$X$ there exists a closed subspace $Y$ of $X$ and a bounded
linear operator from $Y$ to $X$ which is not a compact
perturbation of a multiple of the inclusion map from $Y$ to $X$.
We give sufficient conditions on the spreading models of
seminormalized basic sequences of a Banach space $X$
which guarantee that the ``multiple of the inclusion plus compact''
problem has an affirmative answer for $X$. Our results strengthen
a previous result of the first named author, E.~Odell, Th.~Schlumprecht
and N.~Tomczak-Jaegermann as well as a result of Th.~Schlumprecht.
We give an example of a Hereditarily Indecomposable Banach space where
our results apply. For the proof of our main result we use an extension
of E.~Odell's Schreier unconditionality result for arrays.
\end{abstract}

\bigskip
\setcounter{section}{-1}

\section{Introduction} \label{sectionIntro}

A long-standing open famous question of J.~Lindenstrauss asks whether
on every infinite dimensional Banach space $X$ there exists a (linear
bounded) operator from $X$ to $X$ which is not a compact perturbation
of a multiple of the identity operator on $X$. A weaker question was asked
by T.W.~Gowers in 1996 \cite{G} and by Th.~Schlumprecht in 2003 \cite{S2}:
does every infinite dimensional Banach space $X$ admit a (closed) subspace $Y$
and an operator from $Y$ to $X$ which is not a compact
perturbation of a multiple of the inclusion operator from $Y$ to $X$. We refer
to this question as the ``multiple of the inclusion plus compact'' problem.
The main result of this paper gives sufficient conditions on the spreading models
of seminormalized basic sequences of a Banach space $X$ which ensure that
there exists a subspace $Y$ of $X$ having a basis and an operator from $Y$ to $X$
which is not a compact perturbation of the inclusion map from $Y$ to $X$.

If $X$ and $Y$ are Banach spaces, let ${\mathcal L} (X,Y)$, (respectively ${\mathcal K} (X,Y)$),
denote the set of all (respectively compact) operators from $Y$ to $X$.
If $Y$ is a subspace of $X$ let $i_{Y \to X}$ denote the inclusion map from $Y$ to $X$.
If $Y$ is a subspace of $X$ and $T \in {\mathcal L}(Y,X)$ then the statement
$T \not \in {\mathbb C} i_{Y \to X} + {\mathcal K}(Y,X)$, means that $T$ cannot
be written as a compact perturbation of a multiple of the inclusion map from $Y$ to $X$.
We say that the ``multiple of the inclusion plus compact'' problem has an
affirmative answer on a Banach space $X$ if there exists a subspace $Y$ of $X$
and $T \in \mathcal{L}(Y,X)$ such that $T \not \in {\mathbb C} i_{Y \to X} + {\mathcal K}(Y,X)$.
If $(x_n)_n$ is a basic sequence in a Banach space, let $[(x_n)_n]$
denote the closed linear span of the sequence $(x_n)_n$.

Note that if a Banach space $X$ contains an unconditional basic sequence
$(x_n)_n$ then the operator $T \in {\mathcal L} ([(x_n)_n], X)$ defined by
$T(x_n)= (-1)^n x_n$ does not belong to ${\mathbb C} i_{Y \to X} + {\mathcal K}(Y,X)$.
Thus if a Banach space $X$ contains an unconditional basic sequence then the
``multiple of the inclusion plus compact'' problem has an affirmative answer for $X.$
It is not known whether the above question of Lindenstrauss has an affirmative answer in this case.
Hence for the ``multiple of the inclusion plus compact'' problem we restrict our
attention to Banach spaces with no unconditional basic sequences. Recall that
by the Gowers' dichotomy \cite{G1} every Banach space contains an unconditional
basic sequence or a hereditarily indecomposable (HI) subspace. Recall that a Banach space
$X$ is called HI if no infinite dimensional closed subspace $Y$ of $X$ contains a
complemented subspace $Z$ which is of both infinite dimension and infinite
codimension in $Y$ \cite{GM}. Therefore for the ``multiple of the inclusion plus compact
problem'' we only examine HI saturated Banach spaces.

The ``multiple of the inclusion plus compact'' problem was first studied by Gowers
\cite{G} where he proved that it has an affirmative answer for the HI Banach space
$GM$ which was constructed by Gowers and B.~Maurey \cite{GM}. Moreover, Gowers
conjectured that this problem has an affirmative answer for all reflexive Banach spaces.

Subsequently, Schlumprecht \cite{S2} studied the ``multiple of the inclusion plus compact''
problem and gave sufficient conditions on a Banach space $X$ so that
this problem has an affirmative answer on $X$. One of the main results in
\cite{S2} gives sufficient conditions on the spreading models of weakly null
sequences of an infinite dimensional Banach space $X$ which ensure that the ``multiple
of the inclusion plus compact'' problem has an affirmative answer on $X$. Recall that \cite{BL,BS1,BS2}
for every seminormalized basic sequence $(y_n)$ in a Banach space $X$ and for every
$(\varepsilon _n) \searrow 0$ there exists a subsequence $(x_n)$ of $(y_n)$ and a
seminormalized basic sequence $(\tilde{x}_n)$ (not necessarily in $X$) such that for all
$n \in {\mathbb N}$, scalars $(a_i)_{i=1}^n$ with $|a_i| \leq 1$ and $n \leq k_1 < \cdots < k_n$,
$$
\left| \| \sum_{i=1}^n a_{k_i} x_{k_i} \| - \| \sum_{i=1}^n a_i \tilde{x}_i \| \right| < \varepsilon_n.
$$
The sequence $(\tilde{x}_n)$ is called a spreading model of $(x_n)$. If $X$ is a Banach space then
$SP_w(X)$ will denote the set of spreading models of seminormalized weakly null basic
sequences of $X$. If $(x_n)$ is a seminormalized
weakly null basic sequence then $(\tilde{x}_n)$ is an unconditional basic sequence. Thus if
$X$ is an HI Banach space and $(x_n)$ is a seminormalized basic sequence in $X$ with spreading
model $(\tilde{x}_n)$, then it may be easier to study $(\tilde{x}_n)$ than to study
$(x_n)$ itself. Schlumprecht \cite{S2} introduced the following crucial notion (without assigning a 
name):

\begin{Def} \label{Def:SD}
Let $(x_n)$ and $(z_n)$ be two seminormalized basic sequences (not necessarily in the same Banach
space). For $\varepsilon >0$ define
\begin{equation}
\Delta_{(z_n),(x_n)}(\varepsilon) := \sup \{ \| \sum a_i x_i \| : (a_i) \in c_{00},
| a_i | \leq \varepsilon \text{ and } \| \sum a_i z_i \| \leq 1 \} ,
\end{equation}
where $c_{00}$ denotes the linear space of finitely supported scalar sequences. We say that
$(z_n)$ dominates $(x_n)$ on small coefficients, (denoted by $(x_n) << (z_n)$ and abbreviated as
``$(z_n)$ s.c. dominates $(x_n)$), if
\begin{equation}
\lim_{\varepsilon \searrow 0} \Delta_{(z_n),(x_n)}(\varepsilon) =0.
\end{equation}
\end{Def}

Obviously, if $(z_i)$ and $(x_i)$ are seminormalized basic sequences in some Banach spaces
and $(z_i)>>(x_i)$ then 
\begin{equation} \label{E:A}
 \lim_{\varepsilon \searrow 0}\inf\{\|\sum a_iz_i \|: |a_i|\leq \varepsilon \text{ and }
        \|\sum a_i x_i \| = 1 \} = \infty ,
\end{equation}
where we assume that $\inf \emptyset = \infty$.

Another important notion that was introduced by Schlumprecht \cite{S2} was the following property
which is called ``Property P1'' in the present article.

\begin{Def}
A seminormalized basic sequence $(z_i)$ has Property P1 if
\begin{equation}
\liminf_{n \to \infty} \inf_{A \subset {\mathbb N}, |A| =n} \| \sum_{i \in A} z_i \| =
\infty .
\end{equation}
\end{Def}

One of the main results in \cite{S2} is the following powerful result:

\begin{Thm} \label{Thm:Schlumprecht}
Let $X$ be an infinite dimensional Banach space. Let $(x_n)$ and $(z_n)$ be normalized
weakly null basic sequences in $X$ having spreading models $(\tilde{x}_n)$ and $(\tilde{z}_n)$
respectively, such that $(\tilde{x}_n) << (\tilde{z}_n)$ and $(\tilde{z}_n)$ has Property P1.
Then the ``multiple of the inclusion plus compact'' problem has an affirmative answer on $X$.
\end{Thm}

Another main result in \cite{S2} is its Theorem~$1.4$. The idea of that important
result is that an ordinal index is assigned to every normalized basic sequence of a Banach space,
taking values at most equal to the first uncountable ordinal $\omega_1$. Heuristically speaking,
this index measures how close is the basic sequence that we examine to the unit vector basis of
$\ell_1$. Then an index is assigned to a Banach space as the supremum of the indices of
its normalized weakly null basic sequences. Roughly speaking,
\cite[Theorem~$1.4$]{S2} states that if the index of the
Banach space $X$ is larger than the index of one weakly null normalized basic sequence
in $X$, then the ``multiple of the inclusion plus compact'' problem has an affirmative answer on
$X$.

Another sufficient condition for the ``multiple of the inclusion plus compact'' problem to have
an affirmative solution on a Banach space was given by the first named author,
E.~Odell, Th.~Schlumprecht and N.~Tomczak-Jaegermann \cite{AOST}. This is a sufficient
condition on the spreading models of normalized weakly null basic sequences of a Banach space
$X$. In order to state the mentioned result of \cite{AOST}, we need some more definitions
that we give now. A normalized basic sequence $(x_n)$ is called $1$-subsymmetric if it is
$1$-equivalent to all of its subsequences.

\begin{Def} \label{Krivine}
Let $(x_n)$ be a $1$-subsymmetric basic sequence in some Banach space. The Krivine set
of $(x_n)$ is defined to be the set of all $p$'s in $[1, \infty]$ with the following property.
For all $\varepsilon >0$ and $N \in {\mathbb N}$ there exists $m \in {\mathbb N}$ and scalars
$(\lambda_k)_{k=1}^m$ such that for all scalars $(a_n)_{n=1}^N$,
$$
\frac{1}{1 + \varepsilon} \| (a_n)_{n=1}^N \|_p \leq
\| \sum_{n=1}^N a_n y_n \| \leq (1+ \varepsilon) \| (a_n)_{n=1}^N \|_p
$$
where
$$
y_n = \sum_{k=1}^m \lambda_k x_{(n-1)m+k} \text{ for } n=1, \ldots ,N
$$
and $\| \cdot \|_p$ denotes the norm of the space $\ell_p$.
\end{Def}
The Krivine's theorem as it was proved by H.~Rosenthal \cite{R} and H.~Lemberg \cite{L}
states that if $(x_n)$ is a $1$-subsymmetric basic sequence then the Krivine set of $(x_n)$
is non-empty. In particular, if $(x_n)$ is a seminormalized basic sequence
in a Banach space having spreading model $(\tilde{x}_n)$ then $(\tilde{x}_n)$ is
$1$-subsymmetric, hence the Krivine set of $(\tilde{x}_n)$ is non-empty. The following
result was proved in \cite{AOST}:

\begin{Thm} \label{Thm:AOST}
Let $X$ be a Banach space. Assume that there exist normalized weakly null basic sequences
$(x_n)$, $(z_n)$ in $X$ such that $(x_n)$ has a spreading model $(\tilde{x}_n)$ which is not
equivalent to the unit vector basis of $\ell_1$ and $(z_n)$ has a spreading model
$(\tilde{z}_n)$ such that $1$ belongs to the Krivine set of $(\tilde{z}_n)$. Then the
``multiple of the inclusion plus compact problem'' has an affirmative answer on $X$.
\end{Thm}
This result strengthens the result of Gowers \cite{G} since the Banach space $GM$ satisfies
the condition of Theorem~\ref{Thm:AOST}.

In order to state the main result of our paper, we need to introduce a property which is closely
related to the Property P1 and it is called Property P2 in the present article.
The Property P2 appears without a name in \cite{AOST}.

\begin{Def}
A seminormalized basic sequence $(x_n)$ has Property P2 if
for all $\rho >0$ there exists an $M = M(\rho) \in \N$, such that if $\| \sum a_i z_i \| = 1$
then $|\{i: |a_i| \geq \rho\}| \leq M.$
\end{Def}

One of the main results of our paper is the following:

\begin{Thm}  \label{Thm:main}
Let $X$ be a Banach space containing seminormalized basic sequences $(x_i)_i$ and $(x^n_i)_i$
for all $n \in \N$, such that $0 < \inf_{n,i} \| x^n_i\| \leq \sup_{n,i} \| x^n_i \| < \infty$.
Let  $(z_i)_i$ be a basic sequence not necessarily
in $X$. Assume that $(x_i)$ satisfies:
\begin{equation} \label{eq:condA}
\text{The sequence }(x_i)_i \text{ has a spreading  model } (\tilde{x}_i)_{i=1}^{\infty}
\text{ such that }
               (\tilde{x}_i)_{i\in \N} <<  (z_i)_{i \in \N} .
\end{equation}
Assume that for all $n \in \N$ the sequence $(x^n_i)_i$ satisfies:
\begin{equation} \label{eq:condB}
\begin{split}
& \text{The sequence }(x_i^n)_i \text{ has a spreading model }
                (\tilde{x}_i^n)_i \text{ such that }\\
& (z_i)_{i=1}^n \text{ is $C$-dominated by } (\tilde{x}_i^n)_{i=1}^{n}\text{ for some }C 
\text{ independent of }n.
\end{split}
\end{equation}
Assume that the sequence $(z_i)$ has a spreading model which has Property P2.
Then there exists a subspace $Y$ of $X$ which has a basis and an operator
$T \in \mathcal{L}(Y,X)$ which is not a compact perturbation of a multiple of the inclusion map.
\end{Thm}

Theorem~\ref{Thm:main} obviously implies Theorem~\ref{Thm:Schlumprecht}.
Indeed, if the assumptions of Theorem~\ref{Thm:Schlumprecht} apply for a Banach
space $X$, then let the sequence ``$(z_i)$'' of Theorem~\ref{Thm:main} to be the
sequence $(\tilde{z}_n)$ (which appears in the assumptions of
Theorem~\ref{Thm:Schlumprecht}) and the sequences $(x^n_i)_i$ of
Theorem~\ref{Thm:main} to be all equal to the sequence $(z_n)$ (which appears in the
assumptions of Theorem~\ref{Thm:Schlumprecht}). Notice that since $(z_n)$ is weakly null,
we have that $(\tilde{z}_n)$ is unconditional by \cite{BS1}, \cite{BS2} and it is trivial to verify
that any seminormalized unconditional basic sequence with Property P1 must have Property P2.

Another easy corollary of Theorem~\ref{Thm:main} is obtained if we set $(z_i)$ to be the unit
vector basis of $\ell_p$ for some fixed $p \in [1, \infty)$. Then we obtain the following result
(we present its short proof in Section~\ref{sec:main}).

\begin{Thm} \label{thm:AOST_p}
Let $X$ be a Banach space. Assume that there exist seminormalized
basic sequences $(x_i)$, $(y_i)$ in $X$ such that $(x_i)$ has spreading model $(\tilde x_i)$ which is
s.c. dominated
by the unit vector basis of $\ell_p$, for some $p \in [1,\infty)$ and $(y_i)$ has spreading model
$(\tilde y_i)$ such that $p$ belongs to the
Krivine set of $(\tilde y_i)$. Then there exists a subspace $Y$ of $X$ with a basis
and an operator $T \in \mathcal{L}(Y,X)$ which
is not a compact perturbation of a multiple of the inclusion map.
\end{Thm}

By \cite[Proposition 2.1]{AOST} we know that a seminormalized subsymmetric basic sequence
is not equivalent to the unit vector basis of $\ell_1$ if and only if it is s.c. dominated
by the unit vector basis of $\ell_1$. Thus Theorem~\ref{thm:AOST_p} for $p=1$ implies
Theorem~\ref{Thm:AOST}.

In Section~\ref{sec:set} we present an equivalent statement to the following question.
Given a Banach space $X$ does there exists a subspace $Y$  having a basis
and $T \in \mathcal{L}(Y,X)$ such that $T \not \in \mathbb{C}i_{Y\rightarrow X} + \mathcal{K}(Y,X)$?
All of the above mentioned results in fact assert that this last problem has an affirmative answer
on a Banach space $X$ under the corresponding assumptions on $X$ given by each result.

In Section~\ref{sec:odell} we extend the classical result of Odell on Schreier
unconditionality. Recall
that a finite subset $F$ of $\N$ is called a Schreier set if $|F|\leq \min(F)$ (where $|F|$
denotes the cardinality of $F$).  A basic sequence $(x_n)$
is defined to be {\em Schreier unconditional} if there is a constant $C>0$ such that for all scalars
$(a_i) \in c_{00}$ and for all Schreier sets $F$ we have

$$\| \sum_{i \in F} a_i e_i\| \leq C\| \sum a_i e_i\|.$$

\noindent In this case $(e_i)$ is called {\bf $C$-Schreier unconditional}.
The important notion of Schreier unconditionality was introduced by E.~Odell \cite{O} and has
inspired rich literature on the subject (see for example \cite{AMT}, \cite{DOSZ}). Earlier very similar
results can be found in \cite[page $77$]{MR} and \cite[Theorem $2.1'$]{R3}.

\begin{Thm} \label{thm:odell} \cite{O} 
Let $(x_n)$ be a normalized weakly null sequence in a Banach space.  Then for any
$\varepsilon >0$, $(x_n)$ contains a $(2+\varepsilon)$- Schreier unconditional subsequence.
\end{Thm}

\noindent The main result of Section~\ref{sec:odell} is Theorem~\ref{Th:ASU} where we extend
Theorem~\ref{thm:odell} to arrays.

In Section~\ref{sec:main} we prove Theorem~\ref{Thm:3.1} which also gives sufficient conditions
on a Banach space $X$ so that the ``multiple of the inclusion plus compact'' problem 
has an affirmative
answer on $X$. Then Theorem~\ref{Thm:3.1} is used in the proof of Theorem \ref{Thm:main} which
is further used in the proof of Theorem~\ref{thm:AOST_p}.   The main result
of Section~\ref{sec:odell}, Theorem~\ref{Th:ASU}, plays an important role in the proof of
Theorem~\ref{Thm:3.1}.
In Section~\ref{sec:main} we also examine the relationship between the Properties P1 and P2.

As already mentioned we can restrict the question of the ``multiple of the inclusion plus compact''
problem to HI saturated spaces.  So for a nontrivial application of the above results we need to look
at the
list of HI spaces.  In his 2000 dissertation N.~Dew \cite{D} introduced a new HI space which we
refer to as space $D$. In Section~\ref{sec:dew} we examine some of the basic properties of $D$
and we apply Theorem~\ref{thm:AOST_p} to prove that the``multiple of the inclusion plus compact''
problem has an affirmative answer in $D$.

%------------------------------------------------------------------------------------------------
\section{An almost equivalent reformulation of the ``multiple of the inclusion
        plus compact'' problem } \label{sec:set}
%------------------------------------------------------------------------------------------------

A closely related problem to the ``multiple of the inclusion
plus compact'' problem is the following

\begin{Qst} \label{qst:3}
Assume $X$ is an infinite dimensional Banach space.
Is there a closed subspace $Y$ of $X$ having a basis and an operator $T \in \mathcal{L}(Y,X)$
so that $T \not \in \mathbb{C} i_{Y \rightarrow X} + \mathcal{K}(Y,X)$?\\
\end{Qst}

Notice that Theorems \ref{Thm:main} and \ref{thm:AOST_p} provide sufficient conditions for
this question to have an affirmative answer. Moreover, the proofs of Theorems
\ref{Thm:Schlumprecht} and \ref{Thm:AOST} also reveal that they provide an affirmative answer
to this question.

Before presenting Proposition~\ref{lem:eq} which is the main result of the section, we start
with the following remark which will be used in the proof of Proposition~\ref{lem:eq}.

\begin{Rmk} \label{Lem:uni}
Let $X$ be a Banach space containing no unconditional basic sequence.  Let $(x_n)$
be a seminormalized basic sequence in $X$ and $S \in \mathcal{L}([(x_n)],X)$
such that $(Sx_n)$ converges.  Then there exists a subsequence $(x_{n_k})$
of $(x_{n})$ such that the restriction of $S$ on the span of $(x_{n_{2k}} - x_{n_{2k-1}})$
is compact.
\end{Rmk}

\begin{proof}
Since $X$ does not contain any unconditional basic sequence, no subsequence of $(x_n)$
is equivalent to the unit vector basis of $\ell_1$, hence by Rosenthal's $\ell_1$ Theorem
\cite{R2} there exists a subsequence $(x_{n_k})$ of  $(x_{n})$ which is weak Cauchy.  Thus
the sequence $(x_{n_{2k}} - x_{n_{2k-1}})$ is seminormalized and weakly null.  Since
$(Sx_n)$ converges we have that $(S(x_{n_{2k}} - x_{n_{2k-1}}))$ converges to zero.  By
passing to a further subsequence of $(x_n)$ and relabeling we may assume that
$\sum \|S(x_{n_{2k}} - x_{n_{2k-1}})\| < \infty$ which easily implies that the restriction
of $S$ on the span of $(x_{n_{2k}} - x_{n_{2k-1}})_k$ is compact.  \end{proof}

The next result gives an equivalent characterization
of Question~\ref{qst:3}.

\begin{Prop}  \label{lem:eq}
Let $X$ be a Banach space. The following are equivalent.
\begin{itemize}
\item[(A)] There exists a basic sequence $(x_n)$ in $X$ and an
        operator $T \in \mathcal{L}([(x_n)],X)$ such that
        $T\not \in \mathbb{C} i_{[(x_n)] \rightarrow X} + \mathcal{K}([(x_n)],X)$.

\item[(B)] There exists a seminormalized basic sequence $(x_n)$ in $X$ and a sequence $(y_n)$ in $X$
        such that $(x_n)$ dominates $(y_n)$ and one of the following three happen.
\begin{itemize}
\item[(i)] For all scalars $\lambda$, $(y_n-\lambda x_n)$ has no converging subsequence.

\item[(ii)] There exists a scalar $\lambda$, such that $(y_n-\lambda x_n)$ converges
        and there exists a bounded sequence $(z_n) \subseteq$~span~$(x_i)_{i=1}^\infty$ such that
        $(w_n-\lambda z_n)$ has no converging subsequence, where each
        $w_n \in$~span~$(y_i)_{i=1}^\infty$  has the same distribution with respect to the $(y_i)_i$
        as $z_n$ has with respect to the $(x_i)_i$.

\item[(iii)] There exist scalars $\lambda_1 \not = \lambda_2$ and increasing sequences
        $(k_n^1)$, $(k_n^2)$ of $\N$ such that $(y_{k^i_n} - \lambda_i x_{k^i_n})_n$
        converges for $i \in \{1,2\}$.
\end{itemize}
\end{itemize}
\end{Prop}

\begin{proof}
Note that if the Banach
space $X$ contains an unconditional basic sequence then (A) is satisfied as we noticed in the 
Introduction. Also, in that case, (B)(i) is satisfied (if we set $(x_n)$ to be a seminormalized 
unconditional basic sequence in $X$ and $(y_n)$ to be equal to $((-1)^nx_n)$. Thus we can restrict
our attention to a Banach space $X$ containing no unconditional basic sequence 
(i.e. by Gowers' dichotomy \cite{G1}, $X$ is saturated with HI spaces).  

To show (A) implies (B) let $T \in \mathcal{L}([(x_n)],X)$ for some seminormalized
basic sequence $(x_n) \subseteq X$ where
$T\not \in \mathbb{C} i_{[(x_n)]\rightarrow X} + \mathcal{K}([(x_n)],X)$.  Then
we see that $(x_n)$ dominates $(y_n) := (T(x_n))$.  Then
either $(y_n - \lambda x_n)$ has no converging subsequence for all scalars $\lambda$
(hence (i) holds) or there is a unique scalar $\lambda$ such that $(y_n - \lambda x_n)$
has a converging subsequence or there exist scalars
$\lambda_1 \not = \lambda_2$ and increasing sequences $(k_n^1)$, $(k_n^2)$ of $\N$ such
that $(y_{k^i_n} - \lambda_i x_{k^i_n})_n$ converges for $i \in \{1,2\}$(hence (iii) holds).
Thus if (B) is not valid, then
there exists a unique scalar $\lambda$ such that $(y_n-\lambda x_n)$ converges and for
all bounded sequences $(z_n) \subseteq$ span $(x_n)$ we have $(T(z_n)-\lambda z_n)$
converges.  Thus $T -\lambda i_{[(x_n)] \rightarrow X}$ is compact which is a contradiction.

To show (B) implies (A) we assume we have a pair of sequences $(x_n)$ and $(y_n)$ in $X$ with $(x_n)$ a
seminormalized basic sequence, and $(x_n)$ dominating $(y_n)$.  Define a bounded linear
operator $T:[(x_n)] \rightarrow X$ by $T(x_n) = y_n.$  We show that in each case (i), (ii) and (iii) we
have $T\not \in \mathbb{C} i_{[(x_n)] \rightarrow X} + \mathcal{K}([(x_n)],X)$.

\smallskip

\noindent {\bf Case (i):} Assume $T \in \lambda i_{[(x_n)] \rightarrow X} + \mathcal{K}([(x_n)],X)$ for some $\lambda.$
Then $T-\lambda i_{[(x_n)] \rightarrow X}$ is compact.  Notice that $(x_n)$ is bounded thus we have
$((T-\lambda i_{[(x_n)] \rightarrow X})x_n)_n = (y_n-\lambda x_n)_n$ has a convergent subsequence.  A contradiction to
the assumption (i).

\smallskip

\noindent {\bf Case (ii):}
Assume that for some scalar $\mu$ we have that
$T-\mu i_{[(x_n)] \rightarrow X} \in \mathcal{K}([(x_n)],X)$.
On the other hand, by Remark~\ref{Lem:uni} applied to $S := T-\lambda i_{[(x_n)] \rightarrow X}$
and $(x_n)$ we
have that there exists a subsequence $(x_{n_k})$ of $(x_{n})$ such that the restriction of
$T-\lambda i_{[(x_n)] \rightarrow X}$ on the span of $(x_{n_{2k}}-x_{n_{2k-1}})_k$ is compact.
Let $Y$ denote the span of $(x_{n_{2k}}-x_{n_{2k-1}})_k$.  Then the restriction of the operator

$$(\lambda - \mu) i_{[(x_n)] \rightarrow X} = (T-\mu i_{[(x_n)] \rightarrow X}) -
         (T-\lambda i_{[(x_n)] \rightarrow X})$$

\noindent on $Y$ is compact.  Since $Y$ is infinite dimensional we obtain that $\lambda = \mu$.
Hence $T-\lambda i_{[(x_n)] \rightarrow X} \in \mathcal{K}([(x_n)],X)$.  This contradicts
the assumption (ii) that there exists a bounded sequence $(z_n) \in $ span$(x_i)_i$ such that
$(w_n-\lambda z_n)_n = ((T-\lambda i_{[(x_n)] \rightarrow X})z_n)_n$ has no converging subsequence.

\smallskip

\noindent {\bf Case (iii):}  Assume that for some scalar $\mu$ we have that
$T-\mu i_{[(x_n)] \rightarrow X} \in \mathcal{K}([(x_n)],X)$.  Then as in the
proof of case (ii) we obtain that $\mu = \lambda_1$ and $\mu = \lambda_2$ contradicting
the fact that $\lambda_1 \not = \lambda_2$.
\end{proof}

%---------------------------------------------------------
\section{Extension of Odell's Schreier Unconditionality} \label{sec:odell}

The main result of this section is Theorem~\ref{Th:ASU} which is an extension of
Theorem~\ref{thm:odell} to an array of vectors of a Banach space such that
each row is a seminormalized weakly null sequence.  Then Theorem~\ref{Th:ASU} guarantees
the existence of a subarray which preserves all the rows of the original array and
has a Schreier type of unconditionality.  Theorem~\ref{Th:ASU} will be an important
tool in the proof of the main result of this article (Theorem~\ref{Thm:main}) in 
section~\ref{sec:main}.

We now define the notions of array, subarray and regular array in a Banach space.  An {\bf array}
in a Banach space $X$ is a sequence of vectors in $(x_{i,j})_{i\in \N; j \in J_i}\subseteq X$
where $J_i$ is an infinite subsequence of $\N$ for all $i \in \N$, say
$J_i = \{j_{i,1}<j_{i,2}<\cdots \}$ and $(x_{i,j_{i,k}})_{k \in \N}$ is a seminormalized weakly
null sequence in $X$ for all $i \in \N$.
Let {\bf $<_{r\ell}$} denote the reverse lexicographical
order on $\N^2$.  Let $(x_{i,j})_{i\in \N; j \in J_i}$ be an array in a Banach space $X$.
A {\bf subarray} of $(x_{i,j})_{i\in \N; j \in J_i}$ is an array $(y_{i,\ell})_{i\in \N; \ell \in L_i}$
in $X$ which satisfies the following two properties:

\begin{equation}
 \{y_{i, \ell}: \ell \in L_i \} \subseteq \{x_{i, j}: j \in J_i \} \text{ for all } i \in \N
\end{equation}

\noindent and

\begin{equation}
\begin{split}
 & \text{if }J_i =  \{j_{i,1}<j_{i,2}<\cdots \} \text{, and } L_i =  \{\ell_{i,1}<\ell_{i,2}<\cdots \} \text{ for all }
        i \in \N \text{ then there } \\
 & \text{exists a } <_{r\ell} \text{-order preserving map }
        H: \{(i, j_{i,k}):i,k \in \N\} \rightarrow \{(i, j_{i,k}):i,k \in \N\}\\
 & \text{ such that } y_{i, \ell_{i,k}} = x_{H(i, j_{i,k})} \text{ for all } i,k \in \N.
\end{split}
\end{equation}

\noindent A {\bf regular array} in a Banach space $X$ is an array $(x_{i,j})_{i,j \in \N; i \leq j}$
which is a basic sequence when it is ordered with the reverse lexicographic order: $x_{1,1},x_{1,2},
x_{2,2},x_{1,3},x_{2,3},x_{3,3},x_{1,4},\ldots$.  For convenience, throughout this paper, we denote
the index set of a regular array by $I$, i.e. $I = \{(i,j) \in \N \times \N: i \leq j\}.$ The only reason that
we choose to work with $I$ rather than $\N^2$ is because $I$ has an enumeration (given by the
reverse lexicographic order) that is easy to write down.

Notice the following:

\begin{Rmk} \label{rmk:regular}
If $(x_{i,j})_{(i,j) \in I}$ is a regular array and $(y_{i,j})_{(i,j) \in I}$ is a subarray of $(x_{i,j})_{(i,j) \in I}$
then $(y_{i,j})_{(i,j) \in I}$ is also regular.
\end{Rmk}

The proof of the following remark can be found in functional analysis text books such as
\cite[Theorem 1.5.2]{AK} or \cite[Lemma 1.a.5]{LT}).

\begin{Rmk} \label{rmk:kalton}
Let $(x_i)_{i=1}^N$ be a finite basic sequence in some infinite dimensional Banach space $X$ having basis constant
$C$.  Let $(y_i)$ be a seminormalized weakly null sequence $X$ and $\varepsilon > 0$.  Then there exists an $n \in \N$
such that $(x_1, x_2, \ldots, x_N, y_n)$ is a basic sequence with constant $C(1+\varepsilon)$.
\end{Rmk}

By repeated application of Remark~\ref{rmk:kalton} we obtain the following.

\begin{Rmk} \label{rmk:regular2}
Let $X$ be a Banach space and for every $i \in \N$ let $(x_{i,j})_{j=i}^\infty$ be a seminormalized weakly null
sequence in $X$.  Then there exists a subarray  $(y_{i,j})_{(i,j) \in I}$ of  $(x_{i,j})_{(i,j) \in I}$ which is regular. Moreover, the basis constant of $(y_{i,j})_{(i,j) \in I}$ can be chosen to be 
arbitrarily close to $1$.
\end{Rmk}

For $p \in \N$ any element $\vec a = (a_i)_{i = 1}^p$ of $\R^p$ will
be called a $p$-{\bf{pattern}} and for such $\vec{a}$ define
$|\vec{a}| := p.$  Let $(x_{i,j})_{(i,j) \in I}$ be a regular array in a Banach
space $X$.  Let $f \in X^*$, $k \in \N$, $\vec a = (a_i)_{i=1}^p$ a
$p$-pattern and $F = \{j_1, j_2\ldots,j_p\} \subseteq \N$.  We say
that $f$ has pattern $\vec a$ on $(k,F)$ with respect to
$(x_{i,j})_{(i,j)\in I}$ if $f(x_{k,j_i}) = a_i$ for all $i \in
\{1,2,\ldots , p\}$.

\begin{Lem} \label{lem:A}%------------------------------------------------------------------
Let $(x_{i,j})_{(i,j) \in I}$ be a regular array in a Banach space $X$,
$\vec{a}$ be a $p$-pattern, $\mathcal{F} \subseteq 2Ba(X^*)$, $\delta >0$ and  $i_0, j_0, k_0 \in \N$ with $j_0 \geq i_0$.
Then there exists a subarray $(y_{i,j})_{(i,j) \in I}$ of $(x_{i,j})_{(i,j) \in I}$
such that for any $F \subseteq \{k_0, k_0+1, k_0+2, \ldots \}$, with $(i_0, j_0) <_{r\ell} (k_0, \min(F))$ and
$|F| = p$ we have the following:\\

\noindent If there exists $f \in \mathcal{F}$ having pattern $\vec{a}$ on $(k_0,F)$ with respect
to $(y_{i,j})_{(i,j) \in I}$
then there exists $g \in \mathcal{F}$ having pattern $\vec{a}$ on $(k_0,F)$ with respect
to $(y_{i,j})_{(i,j) \in I}$ and $|g(y_{i_0,j_0})|<\delta$. \\

\noindent Additionally, $(y_{i,j})_{(i,j) \in I}$ can be chosen to satisfy
$y_{i,j} = x_{i,j}$  for all $(i, j) <_{r\ell} (i_0, j_0)$.
\end{Lem}

\begin{proof}
First note that there exists $m \in \N$ such that for all $f \in 2Ba(X^*)$ there exists
$j' \in \{j_0,j_0+1,j_0+2,\ldots,m\}$ with $|f(x_{i_0, j'})| < \delta$.

Otherwise assume that for all $m \in \N$ there exists $x_m^* \in 2Ba(X^*)$ with
$|x_m^*(x_{i_0, j})| \geq \delta$ for $j \in \{j_0,\ldots,m\}$.
By passing to a subsequence and relabeling assume that $(x_m^*)$ converges
weak$^*$ to some $x^* \in 2Ba(X^*)$.
Then $|x^*(x_{i_0, j})|\geq \delta$ for all $j \geq j_0$, which contradicts
that each row, in particular $(x_{i_0, j})_{j={i_0}}^\infty$, is weakly null.

Let N = $\{m+1,m+2,\ldots\}$.  Divide the set $[N]^p$ of all $p$-element subsets of N as follows:
$[N]^p = \cup_{j=j_0}^{m+1} \mathcal{A}_j$ where for $j \in \{j_0, j_0+1,j_0+2,\ldots,m \}$ we set

\begin{equation*}
\begin{split}
\mathcal{A}_j =  \big \{  F \in [N]^p:&\textrm{ there exists }f \in \mathcal{F} \textrm{ having pattern } \vec{a}
                \textrm{ on } (k_0, F) \\
        &\text{ with respect to } (x_{i,j})_{(i,j)\in I}
                \textrm{ and } |f(x_{i_0,j})| < \delta \big \}
\end{split}
\end{equation*}

\noindent and

$$\mathcal{A}_{m+1} = \left \{ F\in [N]^p:\textrm{ there is no }f \in \mathcal{F}\textrm{ having pattern }\vec{a}
        \textrm{ on }(k_0, F) \text{ with respect to } (x_{i,j})_{(i,j)\in I} \right \}. $$

By Ramsey's theorem there exist a subsequence $(m_i)_{i=1}^\infty \in [N]$, and
$j'\in\{j_0,j_0+1,\ldots,m+1\}$ such that $[(m_i)_{i=1}^\infty ]^p \subset \mathcal{A}_{j'}$
(where for an infinite subset $M$ of $\N$, $[M]$ denotes the set of all infinite subsequences of $M$).
We then can pass to a subarray $(y_{i,j})_{(i,j) \in I}$ of $(x_{i,j})_{(i,j) \in I}$ by setting

\begin{equation*}
y_{i,j} = \left \{ \begin{array}{ll}
        x_{i,j} & \text{ if } (i,j) <_{r\ell} (i_0,j_0)\\
        x_{i,j'} & \text{ if }(i,j) = (i_0,j_0) \\
        x_{i,m_j} & \text{ if } (i_0,j_0) <_{r\ell} (i,j).
\end{array}
\right.
\end{equation*}

Then $(y_{i,j})_{(i,j) \in I}$ satisfies the conclusion of the lemma, since if for
some $F \subseteq \{k_0, k_0+1, \ldots\}$ with $(i_0, j_0)<_{r\ell} (k_0, \min(F))$
and $|F| = p$ there exists $f \in \mathcal{F}$ having pattern $\vec{a}$ on $(k_0, F)$
with respect to $(y_{i,j})_{(i,j) \in I}$,
then the integer $j'$ that was obtained by Ramsey's theorem could not be equal to $m+1$,
hence $j' \in\{j_0, j_0+1, \ldots,m\}$ and the definition of $A_{j'}$ gives the conclusion.
\end{proof}

\begin{Lem} \label{lem:B}%------------------------------------------------------------------
Let $(x_{i,j})_{(i,j) \in I}$ a be regular array in a Banach space $X$,
$\vec{A}$ be a finite set of patterns, $\mathcal{F} \subseteq 2Ba(X^*)$, $\delta >0$ and $i_0, k_0 \in \N$.
Then there exists some subarray $(y_{i,j})_{(i,j) \in I}$ of $(x_{i,j})_{(i,j) \in I}$
such that for any $\vec{a}$ in $\vec{A}$, $F \subseteq \{k_0, k_0+1, k_0+2, \ldots \}$,
with $|F| = |\vec{a}|$ and $j_0 \in \N$ with $j_0 \geq i_0$ and $(i_0, j_0) <_{r\ell} (k_0, \min(F))$, we have the following:

\vspace{.125in}
\noindent If there exists $f \in \mathcal{F}$ having pattern $\vec{a}$ on $(k_0,F)$ with respect
to $(y_{i,j})_{(i,j) \in I}$
then there exists $g \in \mathcal{F}$ having pattern $\vec{a}$ on $(k_0,F)$ with respect
to $(y_{i,j})_{(i,j) \in I}$ and $|g(y_{i_0,j_0})|<\delta$.

\vspace{.125in}
\noindent Additionally, we can assume that $x_{i,j} = y_{i,j}$ for all $(i,j)<_{r\ell}(i_0,i_0)$.
\end{Lem}

\begin{proof}
We begin by fixing one particular element $\vec{a}$ in $\vec{A}$.  Now apply Lemma~\ref{lem:A}
for $(x_{i,j})_{(i,j) \in I}$, $\vec{a}$,
$\mathcal{F}$, $\delta$, $i_0$, $k_0$ and $j_0 = i_0$ to obtain some subarray
$(y_{i,j}^{\vec{a},i_0})_{(i,j) \in I}$ of $(x_{i,j}^{\vec{a},i_0})_{(i,j) \in I}$ with the property that
for any $F \subseteq \{k_0, k_0+1, k_0+2, \ldots \}$, with $(i_0, j_0) <_{r\ell} (k_0, \min(F))$ and
$|F| = |\vec{a}|$ we have the following.  If there exists $f \in \mathcal{F}$ having pattern
$\vec{a}$ on $(k_0,F)$ with respect to $(y_{i,j}^{\vec{a},i_0})_{(i,j) \in I}$
then there exists $g \in \mathcal{F}$ having pattern $\vec{a}$ on $(k_0,F)$ with respect
to $(y_{i,j}^{\vec{a},i_0})_{(i,j) \in I}$ and $|g(y_{i_0,j_0}^{\vec{a},i_0})|<\delta$.
Moreover, $y_{i,j}^{\vec{a},i_0}= x_{i,j}$ for all $(i,j)<_{r\ell}(i_0,j_0)$.

We repeat inductively on $j_0$ counting upward from $i_0$.  Thus we next apply
Lemma~\ref{lem:A} to $(y_{i,j}^{\vec{a},i_0})_{(i,j) \in I}$,
$\vec{a}$, $\mathcal{F}$, $\delta$, $i_0$, $k_0$, $j_0 = i_0 + 1$, to obtain
some subarray  $(y_{i,j}^{\vec{a},i_0 + 1})_{(i,j)\in I}$ of $(y_{i,j}^{\vec{a},i_0})_{(i,j)\in I}$ with the property that
for any $F \subseteq \{k_0, k_0+1, k_0+2, \ldots \}$ with $(i_0, j_0) <_{r\ell} (k_0, \min(F))$ and
$|F| = |\vec{a}|$ we have the following.  If there exists $f \in \mathcal{F}$ having pattern
$\vec{a}$ on $(k_0,F)$ with respect to $(y_{i,j}^{\vec{a},i_0+1})_{(i,j)\in I}$
then
\begin{itemize}
\item there exists $g \in \mathcal{F}$ having pattern $\vec{a}$ on $(k_0,F)$ with respect
to $(y_{i,j}^{\vec{a},i_0+1})_{(i,j)\in I}$ and $|g(y_{i_0,i_0}^{\vec{a},i_0 + 1})|<\delta$
(since $y_{i_0,i_0}^{\vec{a},i_0+1}=y_{i_0,i_0}^{\vec{a},i_0}$) and
\item there exists $h \in \mathcal{F}$ having pattern $\vec{a}$ on $(k_0,F)$ with respect
to $(y_{i,j}^{\vec{a},i_0+1})_{(i,j)\in I}$ and $|h(y_{i_0,i_0 + 1}^{\vec{a},i_0 + 1})|<\delta$.
\end{itemize}
\noindent Moreover, $y_{i,j}^{\vec{a},i_0+1}= y_{i,j}^{\vec{a},i_0}$ for all $(i,j)<_{r\ell}(i_0,i_0+1)$.

Continue in this manner for each $j_0 \in \{i_0+2, i_0+3, \ldots\}$.  Note that by fixing the
elements of the subarray for $(i,j) <_{r\ell} (i_0,j_0)$ at each step $j_0$,
there exists a subarray after infinitely many steps which possesses the properties of all the
previous subarrays.  We call this ``limit'' subarray
$(y_{i,j}^{\vec{a}})_{(i,j)\in I}$ and notice it has the property that for any
$F \subseteq \{k_0, k_0+1, k_0+2, \ldots \}$ with $|F| = |\vec{a}|$ and for all $j_0 \in \N$
with $(i_0, j_0) <_{r\ell} (k_0, \min(F))$, we have the following:\\

\noindent If there exists $f \in \mathcal{F}$ having pattern $\vec{a}$ on $(k_0,F)$ with respect
to $(y_{i,j}^{\vec{a}})_{(i,j) \in I}$
then there exists $g \in \mathcal{F}$ having pattern $\vec{a}$ on $(k_0,F)$ with respect
to $(y_{i,j}^{\vec{a}})_{(i,j) \in I}$ and $|g(y_{i_0,j_0}^{\vec{a}})|<\delta$. \\

Notice also that any further subarray of $(y_{i,j}^{\vec{a}})_{(i,j)\in I}$ has this same
property. Then repeat the above process for each $\vec{a} \in \vec{A}$ to obtain the desired array.
\end{proof}

Notice that if $(y_{i,j})_{(i,j) \in I}$ is the result of applying Lemma~\ref{lem:A} to some regular array
and $(z_{i,j})_{(i,j) \in I}$ is a subarray of $(y_{i,j})_{(i,j) \in I}$ then
$(z_{i,j})_{(i,j) \in I}$ does not necessarily retain the property in the conclusion of Lemma~\ref{lem:A}.
However, if $(y_{i,j})_{(i,j) \in I}$ is the result of applying Lemma~\ref{lem:B} to some regular array and
$(z_{i,j})_{(i,j) \in I}$ is a regular subarray of $(y_{i,j})_{(i,j) \in I}$ then $(z_{i,j})_{(i,j) \in I}$
does retain the property in the conclusion of Lemma~\ref{lem:B}.  This idea is summarized in the following
remark.

\begin{Rmk}   \label{rmk:B} %------------------------------------------------------------------
Let $(x_{i,j})_{(i,j) \in I}$ be a regular array in a Banach space $X$,
$\vec{A}$ be a finite set of patterns, $\mathcal{F} \subseteq 2Ba(X^*)$
, $i_0, k_0 \in \N$ and $\delta >0$.  Then there exists a subarray $(y_{i,j})_{(i,j)\in I}$ of
$(x_{i,j})_{(i,j)\in I}$ such that if $(z_{i,j})_{(i,j)\in I}$ is any subarray of $(y_{i,j})_{(i,j)\in I}$,
then for any $\vec{a}$ in $\vec{A}$, $F \subseteq \{k_0, k_0+1, k_0+2, \ldots \}$,
with $|F| = |\vec{a}|$ and $j_0 \in \N$ with $(i_0, j_0) <_{r\ell} (k_0, \min(F))$, we have the following:

\vspace{.125in}
\noindent If there exists $f \in \mathcal{F}$ having pattern $\vec{a}$ on $(k_0,F)$ with respect
to $(z_{i,j})_{(i,j) \in I}$
then there exists $g \in \mathcal{F}$ having pattern $\vec{a}$ on $(k_0,F)$ with respect
to $(z_{i,j})_{(i,j) \in I}$ and $|g(z_{i_0,j_0})|<\delta$.

\vspace{.125in}
\noindent Additionally, we can assume that $x_{i,j} = y_{i,j}$ for all $(i,j)<_{r\ell}(i_0,i_0)$.
\end{Rmk}

\begin{Lem}   \label{lem:C} %------------------------------------------------------------------
Let $(x_{i,j})_{(i,j) \in I}$ be a regular array in a Banach space $X$,
$\vec{A}$ be a finite set of patterns, $\mathcal{F} \subseteq 2Ba(X^*)$
and $\delta >0$.
Then there exists some subarray $(y_{i,j})_{(i,j)\in I}$ of $(x_{i,j})_{(i,j)\in I}$
such that for all  $\vec{a}$ in $\vec{A}$,
$k_0 \in \N$, $F \subseteq \{k_0, k_0+1, k_0+2, \ldots \}$
with $|F| = |\vec{a}|$ and $(i_0,j_0) \in I$ with
$(1,k_0)\leq_{r\ell}(i_0,j_0)<_{r\ell}(k_0,\min(F))$ we have the following:

\vspace{.125in}
\noindent If there exists $f \in \mathcal{F}$ having pattern $\vec{a}$ on $(k_0,F)$ with respect
to $(y_{i,j})_{(i,j)\in I}$
then there exists $g \in \mathcal{F}$ having pattern $\vec{a}$ on $(k_0,F)$ with respect
to $(y_{i,j})_{(i,j)\in I}$ and $|g(y_{i_0,j_0})| < \delta$.
\end{Lem}

\begin{proof}
We will apply Remark~\ref{rmk:B} inductively with the subarray changing at each step,
but $\vec{A}$, $\mathcal{F}$ and $\delta$ remaining as in the hypothesis and $(i_0, k_0)$
cycling through $\N^2$.  We create the final subarray $(y_{i,j})_{(i,j)\in I}$ of
$(x_{i,j})_{(i,j)\in I}$ inductively one column at a time.  At the $j_0$ step of the
induction we create a subarray $(y_{i,j}^{j_0})_{(i,j)\in I}$ of $(y_{i,j+j_0-1}^{j_0-1})_{(i,j)\in I}$
(where for $j_0=1$, $(y_{i,j}^{0})_{(i,j)\in I} = (x_{i,j})_{(i,j)\in I}$ ) and we set
$y_{i,j_0} = y_{i,j_0}^{j_0}$ for $i \in \{1,2,\ldots,j_0\}$.

\noindent {\bf COLUMN 1}:  Apply Remark~\ref{rmk:B} to $(x_{i,j})_{(i,j) \in I}$,
$\vec{A}$, $\mathcal{F}$, $\delta$, $i_0 = 1$, and $k_0 = 1$
to obtain a subarray $(y_{i,j}^1)_{(i,j) \in I}$ with the property that
for all  $\vec{a}$ in $\vec{A}$,
$F \subseteq \{1, 2, \ldots \}$
with $|F| = |\vec{a}|$ such that
$(1,1)<_{r\ell}(k_0,\min(F))$ we have the following:

\vspace{.125in}
\noindent If there exists $f \in \mathcal{F}$ having pattern $\vec{a}$ on $(k_0,F)$ with respect
to $(y_{i,j}^1)_{(i,j)\in I}$
then there exists $g \in \mathcal{F}$ having pattern $\vec{a}$ on $(k_0,F)$ with respect
to $(y_{i,j}^1)_{(i,j)\in I}$ and $|g(y_{1,1}^1)| < \delta$.

\vspace{.125in}
\noindent We then fix column 1 of $(y_{i,j})_{(i,j) \in I}$ by setting $y_{1,1} := y_{1,1}^1$.

\noindent {\bf COLUMN 2}:  Apply Remark~\ref{rmk:B} to
$(y_{i,j+1}^1)_{(i,j) \in I}$, $\vec{A}$, $\mathcal{F}$, $\delta$, successively
for each $(i_0,k_0) \in \{(1,2),(2,1),(2,2) \}$
to obtain a subarray $(y_{i,j}^2)_{(i,j) \in I}$ with the property that
for all  $\vec{a}$ in $\vec{A}$,
$F \subseteq \{2,3, \ldots \}$
with $|F| = |\vec{a}|$, $i_0 \in \{1,2\}$ and $k_0 \in \{1,2\}$ such that
$(i_0,2)<_{r\ell}(k_0,\min(F))$ we have the following:

\vspace{.125in}
\noindent If there exists $f \in \mathcal{F}$ having pattern $\vec{a}$ on $(k_0,F)$ with respect
to $(y_{i,j}^2)_{(i,j)\in I}$
then there exists $g \in \mathcal{F}$ having pattern $\vec{a}$ on $(k_0,F)$ with respect
to $(y_{i,j}^2)_{(i,j)\in I}$ and $|g(y_{i_0,2}^2)| < \delta$.

\vspace{.125in}
\noindent We then fix column 2 of $(y_{i,j})_{(i,j) \in I}$ by setting $y_{i,2} := y_{i,2}^2$
for $i \in \{1,2\}$.

\noindent {\bf COLUMN ${\bf j_0}$}: Apply Remark~\ref{rmk:B} to
$(y_{i,j+j_0-1}^{j_0-1})_{(i,j) \in I}$, $\vec{A}$, $\mathcal{F}$, $\delta$, successively
for each $(i_0,k_0) \in \{(i,j_0): 1 \leq i < j_0 \} \cup \{(j_0,j): 1 \leq j \leq j_0 \}$
to obtain a subarray $(y_{i,j}^{j_0})_{(i,j) \in I}$ with the property that
for all  $\vec{a}$ in $\vec{A}$,
$F \subseteq \{j_0,j_0+1, \ldots \}$
with $|F| = |\vec{a}|$, $i_0 \in \{1,2, \ldots, j_0\}$ and $k_0 \in \{1,2,\ldots, j_0\}$ such that
$(i_0,j_0)<_{r\ell}(k_0,\min(F))$ we have the following:

\vspace{.125in}
\noindent If there exists $f \in \mathcal{F}$ having pattern $\vec{a}$ on $(k_0,F)$ with respect
to $(y_{i,j}^{j_0})_{(i,j)\in I}$ then there exists $g \in \mathcal{F}$ having pattern
$\vec{a}$ on $(k_0,F)$ with respect to $(y_{i,j}^{j_0})_{(i,j)\in I}$ and
$|g(y_{i_0,j_0}^{j_0})| < \delta$.

\vspace{.125in}
\noindent We then fix column $j_0$ of $(y_{i,j})_{(i,j) \in I}$ by setting $y_{i,j_0} := y_{i,j_0}^{j_0}$
for $i \in \{1,2, \ldots, j_0\}$.

Let $\vec{a}$ in $\vec{A}$, $k_0 \in \N$, $F \subseteq \{k_0, k_0+1, k_0+2, \ldots \}$
with $|F| = |\vec{a}|$ and $(i_0,j_0) \in I$ with
$(1,k_0)\leq_{r\ell}(i_0,j_0)<_{r\ell}(k_0,\min(F))$ all be given.
Since $(1,k_0)\leq_{r\ell}(i_0,j_0)$ we have that $k_0 \leq j_0$.  Since $(i_0,j_0) <_{r\ell} (k_0, \min(F))$
we have that there exists a set $G \subseteq \N$ with $|G|=|F|$, $\min(G) \geq k_0$ and
$(y_{k_0,j})_{j \in F} = (y_{k_0,j}^{j_0})_{j \in G}$.  Thus if there exists $f \in \mathcal{F}$
which has pattern $\vec a$ on $(k_0,F)$ with respect to $(y_{i,j})_{(i,j) \in I}$ then $f$
has pattern $\vec a$ on $(k_0,G)$ with respect to $(y_{i,j}^{j_0})_{(i,j) \in I}$, therefore
by the property of $(y_{i,j}^{j_0})_{(i,j) \in I}$ we obtain that there exists $g \in \mathcal{F}$
which has pattern $\vec a$ on $(k_0,G)$ with respect to $(y_{i,j}^{j_0})_{(i,j) \in I}$ and
$|g(y_{i_0,j_0}^{j_0})| < \delta$.  Hence $g$ has pattern $\vec a$ on $(k_0,F)$ with respect
to $(y_{i,j})_{(i,j) \in I}$ and $|g(y_{i_0,j_0})| < \delta$.
\end{proof}

\begin{Lem}  \label{lem:big} %--- Big Lemma  ------
Let $(x_{i,j})_{(i,j) \in I}$ be a regular array in a Banach space $X$,
$\varepsilon >0$, and $k \in \N$.
Then there exists some subarray $(y_{i,j})_{(i,j) \in I}$ of $(x_{i,j})_{(i,j) \in I}$
such that for any pattern $\vec{a}$ in $[-1,1]^p$ for some $p \leq k$, for any $k_0 \in \N$ and
any $F \subseteq \{k_0, k_0 +1, k_0 +2, \ldots \}$ with $|F| = |\vec{a}|$ we have the following:\\

\noindent If there exists $f \in Ba{X^*}$ having pattern $\vec{a}$ on $(k_0,F)$ with respect to $(y_{i,j})_{(i,j) \in I}$
then there exists $g \in (1+\varepsilon)Ba{X^*}$ having pattern $\vec{a}$ on $(k_0,F)$ with respect to $(y_{i,j})_{(i,j) \in I}$
and
$$\sum_{ \{(i,j) \in I: j \geq k_0 \} \backslash \{(k_0, \ell): \ell \in F \}} |g(y_{i,j})|<\varepsilon$$.
\end{Lem}

\begin{proof}  Let $(\delta_j)_{j=0}^\infty \subseteq (0,1)$  such that
\begin{equation} \label{eq:l1sum}
\frac{1}{\inf\limits_{i,j}\| x_{i,j}\|}
\left( 4Ck\delta_0 + \sum_{j=1}^\infty 4Cj\delta_j \right) < \varepsilon
\end{equation}

\noindent where $C$ is the basis constant for the regular array  $(x_{i,j})_{(i,j) \in I}$.
Let $A_0$ be a $\delta_0$ net for $[-1,1]$ containing zero and for each $j \in \N$
choose a $\delta_j$ net $B_j$ for $[-1,1]$ with
$\{0\} \subseteq A_0 \subseteq B_1 \subseteq B_2 \subseteq \cdots$.  Let

\begin{equation}
\vec{A}=\{\vec{a}=(a_i)_{i=1}^p \in A_0^p:\textrm{ where }1\leq~p\leq k\}
\end{equation}
\noindent and
\begin{equation} \label{eq:F}
\mathcal{F} = \{f \in (1+\frac{\varepsilon}{2})Ba(X^*): f(x_{i,j}) \in B_j \text{ for all } (i,j) \in I\},
\end{equation}
\noindent where with out loss of generality we assume $\varepsilon < 2$ so $\mathcal{F} \subseteq 2Ba(X^*)$.
Since $0 \in A_0$ the zero functional is in $\mathcal{F}$ therefore $\mathcal{F}$ is nonempty .

We construct the subarray $(y_{i,j})_{(i,j) \in I}$ of $(x_{i,j})_{(i,j) \in I}$ inductively.
First we will construct a subarray $(y_{i,j}^1)_{(i,j) \in I}$ of $(x_{i,j})_{(i,j) \in I}$,
then for $j \in \N$ for $j \geq 2$ we will construct a subarray $(y_{i,k}^j)_{(i,k) \in I}$
of $(y_{i,k+j-1}^{j-1})_{(i,k) \in I}$.  Once the subarray $(y_{i,k}^j)_{(i,k) \in I}$ has been
constructed we set $y_{i,j} := y_{i,j}^j$ for $1 \leq i \leq j$.  Since $(y_{i,k}^j)_{(i,k) \in I}$
is a subarray of $(y_{i,k+j-1}^{j-1})_{(i,k) \in I}$ and $y_{i,j} = y_{i,j}^j$ for $1\leq i \leq j,$
we have that $(y_{i,j})_{(i,j)\in I}$ is a subarray of $(x_{i,j})_{(i,j)\in I}$.

Apply Lemma~\ref{lem:C} to $(x_{i,j})_{(i,j)\in I}$, $\vec{A}$, $\delta_1$, $\mathcal{F}$ to obtain a subarray
$(y^1_{i,j})_{(i,j)\in I}$ of $(x_{i,j})_{(i,j)\in I}$ such that for all $\vec{a} \in \vec{A}$,
$k_0 \in \N$, $F \subseteq \{k_0, k_0+1, k_0+2, \ldots \}$ with $|F| = |\vec{a}|$ and $(i_0, j_0) \in I$
with $(1,k_0) \leq_{r\ell} (i_0,j_0)<_{r\ell} (k_0, \min(F))$ we have:
if there exists $f \in \mathcal{F}$ having pattern $\vec{a} \in \vec{A}$ on $(k_0,F)$ with respect
to $(y_{i,j}^1)_{(i,j)\in I}$
then there exists $g \in \mathcal{F}$ having pattern $\vec{a}$ on $(k_0,F)$ with respect
to $(y_{i,j}^1)_{(i,j)\in I}$ and $|g(y_{i_0,j_0}^1)| < \delta_1$.
Define the elements of the first column of $(y_{i,j})_{(i,j)\in I}$ by setting
$y_{1,1} := y^1_{1,1}$.  Define for each $b \in B_1$ the set
\begin{equation} \label{eq:Fb}
 \mathcal{F}_b = \{ f \in \mathcal{F}: f(y_{1,1}) = b \}.
\end{equation}
Apply Lemma~\ref{lem:C} to $(y^1_{i,j+1})_{(i,j)\in I}$, $\vec{A}$, $\delta_2$, $\mathcal{F}_b$ successively
for each $b \in B_1$ to obtain a subarray
$(y^2_{i,j})_{(i,j)\in I}$ of $(y^1_{i,j+1})_{(i,j)\in I}$ such that
for all $\vec{a} \in \vec{A}$,
$k_0 \in \N$, $F \subseteq \{k_0, k_0+1, k_0+2, \ldots \}$ with $|F| = |\vec{a}|$ and $(i_0, j_0) \in I$
with $(1,k_0) \leq_{r\ell} (i_0,j_0)<_{r\ell} (k_0, \min(F))$ we have for all $b \in B_1$:
if there exists $f \in \mathcal{F}_b$ having pattern $\vec{a} \in \vec{A}$ on $(k_0,F)$ with respect
to $(y_{i,j}^2)_{(i,j)\in I}$
then there exists $g \in \mathcal{F}_b$ having pattern $\vec{a}$ on $(k_0,F)$ with respect
to $(y_{i,j}^2)_{(i,j)\in I}$ and $|g(y_{i_0,j_0}^2)| < \delta_2  $.

Define the elements of the second column of $(y_{i,j})_{(i,j)\in I}$ by setting
$y_{1,2} := y^2_{1,2}$, and  $y_{2,2} := y^2_{2,2}$.  For
each $\vec{b} = (b_1, b_2, b_3) \in B_1 \times B_2 \times  B_2$ set
\begin{equation} \label{eq:Fbb}
 \mathcal{F}_{\vec{b}} = \{ f \in \mathcal{F}: f(y_{1,1}) = b_1, f(y_{1,2}) =
        b_2\textrm{ and  }f(y_{2,2}) = b_3 \}.
\end{equation}

Apply Lemma~\ref{lem:C} to $(y^2_{i,j+2})_{(i,j)\in I}$, $\vec{A}$, $\delta_3$, $\mathcal{F}_{\vec{b}}$ successively
for each $\vec{b} \in B_1 \times B_2 \times  B_2$ to obtain a subarray
$(y^3_{i,j})_{(i,j)\in I}$ of $(y^2_{i,j+2})_{(i,j)\in I}$ such that
for all $\vec{a} \in \vec{A}$,
$k_0 \in \N$, $F \subseteq \{k_0, k_0+1, k_0+2, \ldots \}$ with $|F| = |\vec{a}|$ and $(i_0, j_0) \in I$
with $(1,k_0) \leq_{r\ell} (i_0,j_0)<_{r\ell} (k_0, \min(F))$ we have for all
$\vec{b} \in  B_{p_{1,1}} \times B_{p_{1,2}} \times  B_{p_{2,2}}$:
if there exists $f \in \mathcal{F}_{\vec{b}}$ having pattern $\vec{a} \in \vec{A}$ on $(k_0,F)$ with respect
to $(y_{i,j}^3)_{(i,j)\in I}$
then there exists $g \in \mathcal{F}_{\vec{b}}$ having pattern $\vec{a}$ on $(k_0,F)$ with respect
to $(y_{i,j}^3)_{(i,j)\in I}$ and $|g(y_{i_0,j_0}^3)| < \delta_3$.
Define the elements in the third column $(y_{i,j})_{(i,j) \in I}$ by setting
$y_{1,3} := y^3_{1,3}$, $y_{2,3} := y^3_{2,3}$and $y_{3,3} := y^3_{2,3}$.
Continue in this manner to create the subarray $(y_{i,j})_{(i,j) \in I}$ of
$(x_{i,j})_{(i,j) \in I}$.

Let $\tilde f \in Ba (X^*)$, $\vec{c}$ be a $p$-pattern for $p \leq k$, $k_0 \in \N$
and $F \subseteq \{k_0,k_0+1,k_0+2,\dots\}$ with $|F| =p$ such that $\tilde{f}$ has pattern
$\vec{c}$ on $(k_0,F)$ with respect to $(y_{i,j})_{(i,j) \in I}$.

First it is easy to see using (\ref{eq:l1sum}) that since $\tilde f \in Ba(X^*)$ there is $f \in \mathcal{F}$ (as defined in (\ref{eq:F}))
such that
\begin{itemize}
\item for all $j \in F$ we have $f(y_{k_0,j}) \in A_0$   and $ |\tilde f(y_{k_0,j}) -f(y_{k_0,j})| \leq \delta_0$, and
\item if we define the finite set $G \subset \N$ by $(y_{k_0,j})_{j \in F}=(x_{k_0,j})_{j \in G}$, then
for all $(i,j) \in I \setminus \{ (k_0,j): j \in G \}$ we have that
$|\tilde{f}(x_{i,j}) - f(x_{i,j})| \leq \delta_j$ and $f(x_{i,j}) \in B_j$.
\end{itemize}
\noindent Let $\vec{a} := (f(y_{k_0,j}))_{j\in F}$ and note that $f$ has pattern $\vec{a}$ on
$(k_0,F)$ with respect to $(y_{i,j})_{(i,j) \in I}$.  We will find a functional
$g \in (1+\frac{\varepsilon}{2})Ba(X^*)$ such that $g$ has pattern
$\vec{a}$ on $F$ with respect to $(y_{i,j})_{(i,j) \in I}$ and
$\sum_{I'\setminus \{(k_0,j):j\in F\} } |g(y_{i,j})| < \varepsilon.$

We proceed to find such functional $g$.  But first a bit of notation,
for a $p$-pattern $\vec{\alpha} = (\alpha_i)_{i=1}^p$
we define its derivative $\vec{\alpha}' = (\alpha_{i-1})_{i=2}^p$.

We will walk through the index set $I'=\{(i,j) \in I: j \geq k_0 \}$
proceeding through this set in $<_{r\ell}$-order
and at each step find a functional $g_{i,j}$ with the property that if $(i,j) \not \in \{(k_0,j):j\in F\}$ then
$|g_{i,j}(y_{i,j})|$ will be small and ``agree'' with the previous functional
on $\{(i', j') \in I': (i',j')<_{r\ell}(i,j)\}$.  If
$(i, j) \in \{(k_0,j): j \in F \}$ then we will not change the previously
defined functional.  We will assume $k_0 \geq 3$ for purposes of demonstrating the
construction, but if $k_0 = 1$ or $2$ then we proceed similarly.

\noindent{\bf STEP ${\mathbf (1,k_0)}$:} Note $(1,k_0) \not \in \{(k_0,j):j\in F\}$ (since $k_0 \geq 3$).
Let
$$
\vec{b}=(f(y_{1,1}), f(y_{2,1}), f(y_{2,2}), f(y_{3,1}), \ldots , , f(y_{k_0-1, k_0-1})).
$$
Then $f \in {\mathcal F}_{\vec{b}}$, $f$ has pattern $\vec{a}$ on $(k_0,F)$ with respect to
$(y_{i,j})_{(i,j) \in I}$, $(y_{i,j})_{(i,j)\in I'}$ is a subarray of $(y^{k_0}_{(i,j)})_{(i,j) \in I'}$
and $(1,k_0) \leq _{r \ell} (1,k_0) <_{r\ell} (k_0, \min (F))$, (the last inequality is valid
since $k_0 \geq 3$). Thus there exists $g_{1,k_0} \in {\mathcal F}_{\vec{b}}$ such that
$g_{1,k_0}$ has pattern $\vec{a}$ on $(k_0,F)$ with respect to $(y_{i,j})_{(i,j) \in I}$ and
$|g_{1,k_0}(y_{1,k_0})|< \delta_{k_0}$. Set $F_{1,k_0} = F$ and $\vec{a}_{1,k_0} = \vec{a}$.

\noindent{\bf STEP ${\mathbf (2,k_0)}$:} Note that
$(2,k_0) \not \in \{(k_0,j):j\in F_{1,k_0}\}$ (since $k_0 \geq 3$).
Let
$$
\vec{b}=(f(y_{1,1}), f(y_{2,1}), f(y_{2,2}), \ldots , f(y_{k_0-1,k_0-1}), g_{1,k_0}(y_{1,k_0})).
$$

\noindent Then $g_{1,k_0} \in {\mathcal F}_{vec{b}}$, $g_{1,k_0}$ has pattern $\vec{a}_{1,k_0}$ on
$(k_0, F_{1,k_0}))$ with respect to $(y_{i,j})_{(i,j) \in I}$ and $(1,k_0) <_{r\ell} (2,k_0) <_{r\ell}(k_0, \min(F))$,
(the last inequality is valid because $k_0 \geq 3$). Thus there exists $g_{2,k_0} \in {\mathcal F}_{\vec{b}}$
such that $g_{2,k_0}$ has pattern $\vec{a}_{1,k_0}$ on $(k_0,F_{1,k_0})$ with respect to
$(y_{i,j})_{(i,j) \in I}$ and $| g_{2,k_0}(y_{2,k_0})|< \delta_{k_0}$.
Set $F_{2,k_0} = F_{1,k_0}$ and $\vec{a}_{2,k_0} = \vec{a}_{1,k_0}$.

We continue similarly until the $(k_0-1,k_0)$ step.  The step $(k_0,k_0)$ is slightly different.  We
separate this step into two different cases depending on whether or not $(k_0,k_0) \in \{(k_0,j):j\in F_{k_0-1,k_0}\}$.

\noindent{\bf STEP ${\mathbf (k_0,k_0)}$:} If $(k_0,k_0) \in \{(k_0,j):j\in F_{k_0-1,k_0}\}$ then set
$g_{k_0,k_0} = g_{k_0-1,k_0}$,
$F_{k_0,k_0} = F_{k_0-1,k_0} \setminus \{k_0\}$ and $\vec{a}_{k_0,k_0} = \vec{a}_{k_0-1,k_0}'$.

If $(k_0,k_0) \not \in \{(k_0,j):j\in F_{k_0-1,k_0}\}$ then $(k_0,k_0)<_{r\ell} (k_0,\min(F))$. Let
$$
\vec{b}=(f(y_{1,1}),f(y_{2,1}), f(y_{2,2}), \ldots , f(y_{k_0-1,k_0-1}),g_{1,k_0}(y_{1,k_0}),
g_{2,k_0}(y_{2,k_0}), \ldots , g_{k_0-1,k_0}(y_{k_0-1,k_0})) .
$$

\noindent Then $g_{k_0-1,k_0} \in {\mathcal F}_{\vec{b}}$, $g_{k_0-1,k_0}$ has pattern $\vec{a}_{k_0-1,k_0}$
on $(k_0,F_{k_0-1,k_0})$ with respect to $(y_{i,j})_{(i,j)\in I}$, $(y_{i,j})_{(i,j)\in I'}$ is a
subarray of $(y^{k_0}_{i,j})_{(i,j)\in I}$ and $(1,k_0) \leq_{r\ell} (k_0,k_0) <_{r\ell} (k_0, \min(F))$.
Thus there exists $g_{k_0,k_0} \in {\mathcal F}_{\vec{b}}$ such that $g_{k_0,k_0}$ has pattern
$\vec{a}_{k_0-1,k_0}$ on $(k_0,F_{k_0-1,k_0})$ with respect to $(y_{i,j})_{(I,j) \in I}$ and
$g_{k_0,k_0}(y_{k_0,k_0}) < \delta_{k_0}$. In this case set $F_{k_0,k_0} = F_{k_0-1,k_0}$ and
$\vec{a}_{k_0,k_0} = \vec{a}_{k_0-1,k_0}$.

Then start again with the first entry $(1,k_0+1)$ of the next column as in steps $(1,k_0)$ and $(2,k_0)$.

\noindent{\bf STEP $ {\mathbf(1, k_0+1)}$:} Note that $(1,k_0+1) \not \in \{(k_0,j):j\in F_{k_0,k_0}\}$,
(since $k_0 \geq 3$). Let
$$
\vec{b}= ( f(y_{1,1}), f(y_{2,1}), f(y_{2,2}), \ldots , g_{1,k_0}(y_{1,k_0}), g_{2,k_0}(y_{2,k_0}),
g_{k_0,k_0}(y_{k_0,k_0})).
$$

\noindent Then $g_{k_0,k_0} \in {\mathcal F}_{\vec{b}}$, $g_{k_0,k_0}$ has pattern $\vec{a}_{k_0,k_0}$
on $(k_0, F_{k_0,k_0})$ with respect to $(y_{i,j})_{(i,j) \in I}$, $(y_{i,j})_{(i,j) \in I'}$ is a
subarray of $(y^{k_0}_{i,j})_{(i,j) \in I}$ and $(1,k_0) \leq _{r\ell} (1,k_0+1) <_{r\ell} (k_0, \min(F))$,
(the last inequality is valid since $k_0 \geq3$). Thus there exists $g_{1,k_0+1}  \in {\mathcal F}_{\vec{b}}$
such that $g_{1,k_0+1}$ has pattern $\vec{a}_{k_0,k_0}$ on $(k_0,F_{k_0,k_0})$ with respect to
$(y_{i,j})_{(i,j) \in I}$ and $|g_{1,k_0+1}(y_{1,k_0+1})| < \delta_{k_0+1}$.
Set $F_{1,k_0+1} = F_{k_0,k_0}$ and $\vec{a}_{1,k_0+1} = \vec{a}_{k_0,k_0}$.

Continue in this manner to generate a sequence of functionals $(g_{i,j})_{(i,j)\in I'}$.  We only
need to distinguish two cases every time we reach the $k_0$ row as in step $(k_0,k_0)$. Let
$g \in (1+\frac{\varepsilon}{2})Ba(X^*)$ be a weak$^*$-accumulation point of sequence $(g_{i,j})_{(i,j)\in I'}$.
Note $g$ has the following two properties:

\begin{itemize}
\item $g$ has pattern $\vec{a}$ on $(k_0,F)$ with respect to $(y_{i,j})_{(i,j)\in I}$, and
\item $\sum_{(i,j) \in I' \backslash \{(k_0, \ell): \ell \in F \}} |g(y_{i,j})|<\varepsilon$.
\end{itemize}

\noindent Since $g(y_{k_0,\ell}) = f(y_{k_0,\ell}) \in A_0$ for all $\ell \in F$ and
$|f(y_{k_0,\ell})-  \tilde f(y_{k_0,\ell})| < \delta_0$ for all $\ell \in F$,
(\ref{eq:l1sum}) implies that there exists $\tilde g \in X^*$ such that
$\|\tilde g\| \leq \|g\| + \frac{\varepsilon}{2} \leq 1 + \varepsilon$ and
$\tilde g(y_{k_0,\ell}) = \tilde f (y_{k_0,\ell}) \in A_0$ for all $\ell \in F$
(thus $\tilde g$ has pattern $\vec{c}$ on $(k_0,F)$ with respect to $(y_{i,j})_{(i,j) \in I}$)
and

$$  \sum_{(i,j)\in I' \setminus \{(k_0,\ell): \ell \in F\}} |\tilde g (y_{i,j})| < \varepsilon $$

\noindent completing the proof.
\end{proof}

%--------------------------------------------------------------
\begin{Thm} \label{Th:ASU}
Let $(x_{i,j})_{(i,j) \in I}$ be a regular array in a Banach space $X$,
$(M_j)_{j \in \N} \subseteq \N$ be an increasing sequence of integers and $\varepsilon > 0$.
Then there exists a regular subarray $(y_{i,j})_{(i,j) \in I}$ of $(x_{i,j})_{(i,j) \in I}$
such that for any finitely supported scalars $(a_{i,j})_{(i,j) \in I}$, $k_0 \in \N$ and
$F \subseteq \N$ with $|F| \leq M_{\min(F)}$ and $k_0 \leq \min(F)$
we have
$$ \| \sum_{(i,j) \in I} a_{i,j}y_{i,j}  \| \geq
        \frac{1}{2+\varepsilon} \| \sum_{j \in F} a_{k_0,j}y_{k_0,j}  \|. $$
\end{Thm}
\begin{proof}
Let $\eta >0$  such that
\begin{equation} \label{eq:eta}
2(2\eta +1) \leq 2+ \frac{\varepsilon}{2}
\end{equation}
where $C$ is the basis constant of the regular array $(x_{i,j})_{(i,j) \in I}$.
Apply Lemma~\ref{lem:big} to $(x_{i,j})_{(i,j) \in I}$, $\eta$ and $M_1$
to get $(y_{i,j}^1)_{(i,j) \in I}$.  Define $y_{1,1} := y_{1,1}^1$.

Apply Lemma~\ref{lem:big} to $(y_{i,j+1}^1)_{(i,j) \in I}$, $\eta$ and $M_2$
to get $(y_{i,j}^2)_{(i,j) \in I}$.  Define $y_{i,2} := y_{i,2}^2$ for $i = 1,2$.

Assuming that $(y_{i,j}^{\ell-1})_{(i,j) \in I}$ has been defined (and thus $(y_{i,j})_{(i,j)\in I; j < \ell})$
has also been defined) apply Lemma~\ref{lem:big} to $(y_{i,j+\ell -1}^{\ell-1})_{(i,j) \in I}$, $\eta$ and $M_\ell$
to get $(y_{i,j}^{\ell})_{(i,j) \in I}$.  Define $y_{i,\ell} := y_{i,\ell}^\ell$ for $i = 1,2,\ldots,\ell$.

Inductively construct the entire array $(y_{i,j})_{(i,j) \in I}$ and notice $(y_{i,j})_{(i,j) \in I}$
is regular by Remark~\ref{rmk:regular}.

Let $k_0 \in \N$, $F \subseteq \{k_0, k_0+1, k_0+2,\ldots \}$  with $|F| \leq M_{\min(F)}$
and finitely supported scalars $(a_{i,j})_{(i,j) \in I}$ be given.  We can assume without loss of generality that

$$\|\sum_{(i,j) \in I} a_{i,j}y_{i,j} \|  = 1.$$

Then $|a_{i,j}| \leq 2C$ for $(i,j) \in I$.  Let $f \in Ba(X^*)$ such that

$$ f \left ( \sum_{j \in F} a_{k_0,j}y_{k_0,j}  \right )
 = \|\sum_{j \in F} a_{k_0,j}y_{k_0,j} \|.$$

\noindent Let $\vec{a} = (f(y_{k_0,j}))_{j \in F}$ be a $p$-pattern where $p=|F|.$
Obviously $f$ has pattern $\vec{a}$ on $(k_0,F)$ with respect to
$(y_{i,j})_{(i,j) \in I}$.  Then by considering the subarray $(y_{i,j})_{(i,j) \in I, j \geq \min(F)}$
of $(y^{\min(F)}_{i,j})_{(i,j) \in I:j \geq \min(F)}$ we obtain by the above that
there exists $g \in (1+\eta)Ba(X^*)$
having pattern $\vec{a}$ on $(k_0,F)$ with respect to $(y_{i,j})_{(i,j) \in I}$
and $\sum_{ \{(i,j) \in I : j \geq \min(F) \}\backslash  \{(k_0,\ell): \ell\in F\} } |g(y_{i,j})|<\eta$.
Thus

\begin{equation} \label{eq:est}
\begin{split}
1 &= \| \sum a_{i,j}y_{i,j}\| \geq \frac{1}{2C} \| \sum_{j \geq \min(F)} a_{i,j}y_{i,j}\|
                \geq  \frac{1}{2C(1+\eta)} \bigg | \sum_{j \geq \min(F)} a_{i,j} g(y_{i,j}) \bigg | \\
  &\geq  \frac{1}{2C(1+\eta)} \bigg | \sum_{j \in F} a_{k_0,j} g(y_{k_0,j}) \bigg |
                - \frac{1}{2C(1+\eta)} \bigg | \sum_{\{(i,j) \in I: j \geq \min(F) \} \setminus \{(k_0,j):j \in F \} } a_{i,j} g(y_{i,j}) \bigg | \\
  &\geq  \frac{1}{2C(1+\eta)} \| \sum_{j \in F} a_{k_0,j}y_{k_0,j} \|
                - \frac{1}{2C(1+\eta)} \sum_{\{(i,j) \in I: j \geq \min(F) \} \setminus \{(k_0,j):j \in F \} } |a_{i,j}| |g(y_{i,j})| \\
  &\geq  \frac{1}{2C(1+\eta)} \| \sum_{j \in F} a_{k_0,j}y_{k_0,j} \|
                - \frac{1}{1+\eta} \sum_{\{(i,j): j \geq \min(F) \} \setminus \{(k_0,j):j \in F\} } |g(y_{i,j})| \\
  &\geq  \frac{1}{2C(1+\eta)} \| \sum_{j \in F} a_{k_0,j}y_{k_0,j} \|
                - \frac{1}{1+\eta} \eta .
\end{split}
\end{equation}

Thus by (\ref{eq:eta}) and (\ref{eq:est}) we have
$$\|\sum_{(i,j) \in I} a_{i,j}y_{i,j} \| = 1 \geq \frac{1}{2C(2\eta + 1)}
\|\sum_{j \in F} a_{k_0,j}y_{k_0,j}\| \geq \frac{1}{C(2 + \frac{\varepsilon}{2})}\|\sum_{j \in F} a_{k_0,j}y_{k_0,j}\|$$.

Since we can choose $C$, the basis constant of our regular array, arbitrarily close to $1$ 
(see Remark~\ref{rmk:regular2}) we have shown the result.
\end{proof}

%---------------------------------------------------------

\bigskip
\section{Existence of Non-trivial Operators} \label{sec:main}

In this section we will prove Theorem~\ref{Thm:main} which is one of the main results of the paper.
Theorem~\ref{Th:ASU} will play an important role in its proof
(see the proof of Lemma~\ref{lem:6.3}).

For the proof of  Theorem~\ref{Thm:main} we will need the following result which also gives sufficient
conditions for a Banach space $X$ so that the ``multiple of the inclusion plus compact'' 
problem to have an affirmative answer in $X$.

\begin{Thm} \label{Thm:3.1}
Let $X$ be a Banach space containing seminormalized basic sequences $(\underline{x}_i)_i$ and $(x^n_i)_i$
for all $n \in \N$, such that $0 < \inf_{n,i} \| x^n_i\| \leq \sup_{n,i} \| x^n_i \| < \infty$.
Let  $(x_i)$ and $(z_i)_i$ be seminormalized basic sequences not necessarily
in $X$. Assume the following:
\begin{itemize}
\item The sequence $(\underline{x}_i)$ is dominated by the sequence $(x_i)$.
\item The sequence $(x_i)$ satisfies condition (\ref{eq:condA}) of Theorem~\ref{Thm:main}
and has Property P2.
\item For all $n \in \N$ the sequence $(x^n_i)_i$ satisfies condition (\ref{eq:condB})
of Theorem~\ref{Thm:main}.
\end{itemize}
Then there exists a subspace $Y$ of $X$ which has a basis and an operator
$T \in \mathcal{L}(Y,X)$ which is not a compact perturbation of a multiple of the inclusion map.
\end{Thm}

In order to prove Theorem~\ref{Thm:3.1} we need the following two lemmas whose proofs are postponed.

\begin{Lem}  \label{lem:6.2}  %--------------------------------------------------------------
Let $X$ be a Banach space, $(x_n)$ be a seminormalized basic sequence in $X$ having 
Property P2 and
$(z_n)$ be a seminormalized basic sequence not necessarily in $X$. Assume that the sequence $(x_n)$
satisfies condition~(\ref{eq:condA}) of Theorem~\ref{Thm:main}.
Then for all $(\delta_n)_{n=2}^\infty \subseteq (0,\infty)$ there exists an increasing
sequence $M_1<M_2<\cdots$ of positive integers and a subsequence $(x_{n_i})$ of $(x_{i})$ such that
for all $(a_i) \in c_{00}$,
\begin{equation} \label{eq:6.2}
 \|\sum a_i x_{n_i}\| \leq \sup_{n \in \N} \sup_{n \leq F \subseteq \N; |F|\leq M_n} \delta_n \|\sum_{i \in F} a_i z_i\|,
\end{equation}

\noindent for some $\delta_1$ (where ``$n \leq F$'' means $n \leq \min(F)$).
\end{Lem}

\begin{Lem}  \label{lem:6.3}  %--------------------------------------------------------------
Let $X$ be a Banach space, $(\delta_n)$ be a summable sequence of positive numbers,
$(M_n)_{n=1}^\infty \subseteq \N$ be a sequence of positive integers,
$(z_n)$ be a seminormalized basic sequence (not necessarily in $X$)
and for every $n \in \N$ let $(x_j^n)_{j=1}^\infty$ be a weakly null basic sequence in $X$
having spreading model $(\tilde x_j^n)_{j=1}^\infty$ such that
$0< \inf_{n,j}\|x_j^n\| \leq \sup_{n,j} \|x_j^n\| < \infty$ and
condition~(\ref{eq:condB}) of Theorem~\ref{Thm:main} is satisfied.
Then there exists a seminormalized weakly null basic sequence $(y_{i})$ in $X$ such that

\begin{equation} \label{eq:6.3}
 \frac{1}{6C}\sup_n \sup_{n \leq F \subseteq \N; |F|\leq M_n} \delta_n
        \|\sum_{i \in F} a_i z_i\| \leq \|\sum a_i y_i\|.
\end{equation}

\noindent Moreover, $\| y_j\| \geq \frac{\delta_1}{2} \inf_{n,m} \| x^n_m \|$. Furthermore,
if $(x_i^*)$ is any given sequence of functionals in $X^*$ and $\varepsilon > 0$
we can choose $(y_i)$ to satisfy $|x_i^*y_i| < \varepsilon$.
\end{Lem}

We now present the

\begin{proof}[Proof of Theorem~\ref{Thm:3.1}]
Note that the assumptions of 
Lemmas~\ref{lem:6.2} and \ref{lem:6.3} are almost included in the assumptions of 
Theorem~\ref{Thm:3.1}, with the exception that in Lemma~\ref{lem:6.3}
the sequence $(x^n_i)_i$ is assumed to be weakly null for all $n$. We will replace the sequences
$(\underline{x}_i)$, $(x^n_i)_i$, $(x_i)$ and $(z_i)$ by sequences $(\underline{X}_i)$, $(X^n_i)$,
$(X_i)$  and $(Z_i)$ respectively, such that $(X_i)_i$, $(X^n_i)_i$ and $(Z_i)_i$
satisfy the assumptions of Lemma~\ref{lem:6.2} and \ref{lem:6.3}. Notice that if $\ell_1$ embeds
in $X$ then $X$ contains an unconditional basic sequence thus as we mentioned in the Introduction,
the conclusion of Theorem~\ref{Thm:3.1} is valid in this case. Therefore we can assume that 
$\ell_1$ does
not embed in $X$. Then by Rosenthal's $\ell_1$ Theorem \cite{R2} and a diagonal argument, by passing
to subsequences of $(x_i)$, $(x^n_i)_i$ and $(z_i)$ and relabeling, we can assume that 
$(\underline{x}_i)_i$ and
$(x^n_i)_i$ are weakly Cauchy for all $n \in \N$. Let 
$(\underline{X}_i):=(\underline{x}_{2i}-\underline{x}_{2i-1})$, $(X^n_i):=(x^n_{2i}-x^n_{2i-1})$,
$(X_i)_i:=(x_{2i}-x_{2i-1})$
and $(Z_i):=(z_{2i}-z_{2i-1})$. It is trivial to check that all the assumptions of Lemma~\ref{lem:6.2}
and \ref{lem:6.3} are satisfied for the sequences $(X_i)$, $(X^n_i)_i$ and $(Z_i)$. Thus assume that
this is the case for the original sequences $(x_i)$, $(x^n_i)_i$ and $(z_i)$.
Let $(\delta_n)_{n=2}^\infty$ be a
summable sequence of positive numbers.  First apply Lemma~\ref{lem:6.2} to obtain a subsequences
$(x_{n_i})$, $\delta_1 > 0$ and an increasing sequence $(M_n)_{n \in \N}$
of positive integers which satisfies (\ref{eq:6.2}).  For every $i \in \N$ let a norm $1$ functional
$x_i^*$ satisfying $x_i^*\underline{x}_{n_i} = \| \underline{x}_{n_i}\|$.
Then apply Lemma~\ref{lem:6.3} for $(\delta_n)_{n \in \N}$ and $(M_n)_{n \in \N}$
to obtain a basic sequence $(y_i)$ which satisfies (\ref{eq:6.3}).

Assume also that $(y_i)$ satisfies the ``furthermore'' part of the statement of Lemma~\ref{lem:6.3}
for the sequence $(x_i^*)$ and 
$\varepsilon = \delta_1 (\inf_{n,i}\|x_i^n\|)^2/(8 \sup_{n,i}\| x^n_i \|)$.
Note that if
$|\lambda| \geq \frac{4 \sup_{n,i}\| x^n_i \|}{\delta_1 \inf_{n,i}\|x_i^n\|}$ then
$$\|\underline{x}_{n_i} + \lambda y_i\| \geq | \lambda |\|y_i\| - \|\underline{x}_{n_i}\|  \geq
               \frac{4 \sup_{n,i}\| x^n_i \|}{\delta_1 \inf_{n,i}\|x_i^n\|} \frac{\delta_1}{2} \inf_{n,i}\|x^n_i\|
               - \| \underline{x}_{n_i}\|  \geq \| x^n_i\| \geq \inf_{n,i} \| x^n_i\|  $$
\noindent (by the ``moreover'' part of the statement of Lemma~\ref{lem:6.3}).  
Also if $|\lambda| < \frac{4 \sup_{n,i}\| x^n_i \|}{\delta_1 \inf_{n,i}\|x_i^n\|}$
then
$$\|\underline{x}_{n_i} + \lambda y_i\| \geq |x_i^*(\underline{x}_{n_i} + \lambda y_i)| \geq
                \| \underline{x}_{n_i}\| - \frac{4 \sup_{n,i}\| x^n_i \|}{\delta_1 \inf_{n,i}\|x_i^n\|}\varepsilon
                \geq \frac{1}{2} \inf_{n,i} \| x^n_i\| . $$
\noindent Thus for all scalars $\lambda$ we have

$$ \| \underline{x}_{n_i} +\lambda y_{i}  \| \geq \frac{1}{2} \inf_{n,i}\| x^n_i\|.$$

\noindent Thus if we define $T:[(y_{i})] \rightarrow X$ by
$$ T(\sum a_i y_{i}) = \sum a_i \underline{x}_{n_i}$$
\noindent we have that this operator is bounded by (\ref{eq:6.2}), (\ref{eq:6.3})
and our assumption that $(\underline{x}_i)_i$ is dominated by $(x_i)_i$.
We also have that for any scalar $\lambda$,
$(T-\lambda i_{[(y_i)] \to X})(y_k) = \underline{x}_{n_k}-\lambda y_{k}$.  But since $(y_k)$
is weakly null and $x_{n_k}-\lambda y_{k}$  is not norm null, $T-\lambda i_{[(y_i)] \to X}$ 
is not compact.  In other
words $T$ is not a compact perturbation of a scalar multiple of the inclusion.
\end{proof}

Now we present the proof of Lemma~\ref{lem:6.2}.
A less general version of this  lemma can be found in \cite[Lemma 2.4 (a) $\Rightarrow$ (d)]{S2}.
Schlumprecht  assumes that the basic sequence $(z_i)$ is subsymmetric and satisfies Property P1
and we  assume that the
sequence $(x_i)$ has Property P2, which in view of Proposition~\ref{P:P1usm} can be replaced by the
assumption that $(z_i)$ has a spreading model which is unconditional and satisfies Property P1. 
Also Schlumprecht shows the
result for some sequence $(\delta_n)$ while we show it for an arbitrary $(\delta_n)$.
Additionally, we use  different techniques than the ones used in \cite{S2}. Our arguments resemble
the ones found in \cite{AOST}.

%--------------------------------------Proof of Lemma~\ref{lem:6.2}-------------------------------
\begin{proof}[Proof of Lemma~\ref{lem:6.2}]
Since $(x_n)$ has Property P2, for each $\rho >0$ we can define
$M= M(\rho)$ such that if $\|\sum a_ix_i\| = 1 $ then $|\{i:|a_i|>\rho\}|\leq M$.

Let $(\varepsilon_j)_{j=1}^\infty$ be such that

$$ \sum_{j=2}^\infty \frac{\varepsilon_{j-1}}{\delta_j} \leq \frac{1}{2}.$$

Since $(z_n) >> (\tilde x_n)$ by (\ref{E:A}) we may choose
a decreasing sequence $(\rho_j)_{j=1}^\infty \subseteq (0,1]$ such that
$\sum_{j} \sqrt{\rho_j}(j+1) \leq 1/4$ and satisfying the following: for all $(a_i) \in c_{00}$
with $|a_i|\in [0,\sqrt{\rho_j}]$ for each $i$ and $\|\sum a_i \tilde x_i\| = 1$
we have
\begin{equation} \label{eq:star}
\|\sum a_i \tilde x_i\| \leq \varepsilon_j\|\sum a_i z_i\|.
\end{equation}
\noindent Finally let $M_j = M(\rho_j)$ as above.

By the definition of spreading models, by passing to a subsequence of $(x_i)$ and relabeling, we can assume
that if $j\leq F$ and $|F|\leq M_j$ then for all $(a_i) \in c_{00}$,

\begin{equation}\label{eq:eqtosp}
\frac{1}{2}\|\sum_{i \in F} a_i x_i\| \leq \|\sum_{i \in F} a_i \tilde x_i\|
        \leq 2 \|\sum_{i \in F} a_i x_i\|.
\end{equation}

Now fix $(a_i) \in c_{00}$ such that $\| \sum a_ix_i \| = 1$.  For $j \in  \N$ consider the vector
$\tilde y = \sum_{i>j; \rho_j < |a_i| \leq \rho_{j-1}}   a_i \tilde x_i$.  If
$\|\tilde y\| \geq \sqrt{\rho_{j-1}}$ then

\begin{equation*}
\begin{split}
\|\tilde y\| =& \|\tilde y\| \|\frac{\tilde y}{\|\tilde y\| }\|
        = \|\tilde y\| \bigg \| \mathop{\sum_{i>j}}_{\rho_j < |a_i| \leq \rho_{j-1}}
        \frac{a_i}{\|\tilde y\|} \tilde x_i \bigg\| \\
  \leq& \|\tilde y\| \varepsilon_{j-1} \bigg \| \mathop{\sum_{i>j}}_{\rho_j < |a_i| \leq \rho_{j-1}}
        \frac{a_i}{\|\tilde y\|} z_i \bigg \|
        \text{\phantom{a}(by (\ref{eq:star}) since } \bigg \| \mathop{\sum_{i>j}}_{\rho_j < |a_i| \leq \rho_{j-1}}
        \frac{a_i}{\|\tilde y\|} \tilde x_i \bigg \| = 1)\\
  =& \varepsilon_{j-1} \bigg \| \mathop{\sum_{i>j}}_{\rho_j < |a_i| \leq \rho_{j-1}}  a_i z_i  \bigg \|.
\end{split}
\end{equation*}

\noindent Thus in general, (without assuming that $\|\tilde y\|\geq\sqrt{\rho_{j-1}}$), we get
\begin{equation} \label{eq:rho}
\|\tilde y\| \leq \sqrt{\rho_{j-1}} +
\varepsilon_{j-1} \|\sum_{i>j; \rho_j < |a_i| \leq \rho_{j-1}}  a_i z_i \|.
\end{equation}

Let $\rho_0$  be twice the basis constant of $(x_i)$ divided by the $\inf \|x_i\|$.
Since $\|\sum a_ix_i\|=1$, we have that $|a_i| \leq \rho_0$.
\begin{equation*}
\begin{split}
1 =& \| \sum a_ix_i \| \leq      \sum_{j=1}^\infty \|\sum_{\rho_j < |a_i| \leq \rho_{j-1}} a_ix_i \| \\
  \leq& \|\sum_{\rho_1 < |a_i| \leq \rho_0} a_ix_i \| +
                \sum_{j=2}^\infty \|\sum_{i\leq j; \rho_j < |a_i| \leq \rho_{j-1}} a_ix_i \| +
                \sum_{j=2}^\infty \|\sum_{i > j; \rho_j < |a_i| \leq \rho_{j-1}} a_ix_i \| \\
  \leq& \sup_{F\subseteq \N, |F|\leq M_1} \delta_1 \|\sum_{i\in F} a_iz_i\|+
                \sum_{j=2}^\infty j \rho_{j-1} +
                2\sum_{j=2}^\infty \|\sum_{i > j; \rho_j < |a_i| \leq \rho_{j-1}} a_i \tilde x_i \|
\end{split}
\end{equation*}
\noindent where
$$
\delta_1 = \sup \left\{ \frac{\|\sum_{i \in F} a_i x_i\|}{\|\sum_{i \in F} a_i z_i\|}
:(a_i)_{i \in F} \subseteq {\mathbb C} \text{ with } |F| \leq M_1 \text{ and }(a_i)_{i \in F} \not = 0^F \right\}$$
which is clearly finite by using an $\ell_1$  estimate for the numerator and an  $\ell_\infty$ estimate
for the denominator. Note the third piece of the last inequality is true by (\ref{eq:eqtosp}) since
the cardinality of $\{i > j: \rho_j < |a_i| \leq \rho_{j-1}\}$ is at
most $M_j$. Continuing the calculations from above, we get
\begin{equation*}
\begin{split}
1  \leq& \sup_{F\subseteq \N, |F|\leq M_1} \delta_1\|\sum_{i\in F} a_iz_i\|+
                \frac{1}{4} +
                2\sum_{j=2}^\infty \sqrt{\rho_{j-1}} +
                2\sum_{j=2}^\infty  \varepsilon_{j-1} \|\sum_{i>j; \rho_j < |a_i| \leq \rho_{j-1}}  a_i z_i \|
                \text{ by (\ref{eq:rho})}\\
  \leq& \sup_{F\subseteq \N, |F|\leq M_1} \delta_1\|\sum_{i\in F} a_iz_i\|+
                \frac{1}{2} +
                2\sum_{j=2}^\infty  \frac{\varepsilon_{j-1}}{\delta_j}\delta_j
                        \|\sum_{i>j; \rho_j < |a_i| \leq \rho_{j-1}}  a_i z_i \| \\
  \leq& \sup_{F\subseteq \N, |F|\leq M_1} \delta_1\|\sum_{i\in F} a_iz_i\|+
                \frac{1}{2} +
                \sup_{n \geq 2} \sup_{n \leq F \subseteq \N; |F|\leq M_n} \delta_n \|\sum_{i \in F} a_i z_i\|.
\end{split}
\end{equation*}

\noindent Thus
$$ 1 \leq 2\sup_{n\in \N}
        \sup_{n \leq F \subseteq \N; |F|\leq M_n} \delta_n \|\sum_{i \in F} a_i z_i\|\\$$
\noindent proving the lemma.
\end{proof}

Now we present the proof of Lemma~\ref{lem:6.3}.

%-------------------------------Proof of Lemma~\ref{lem:6.3}----------------------------------
\begin{proof}[Proof of Lemma~\ref{lem:6.3}]
Assume that for each $n \in \N$ $(x_j^{M_n})_j$ has spreading model $(\tilde x_j^{M_n})_{j=1}^{M_n}$,
$(\tilde x_j^{M_n})_j$ $C$-dominates $(z_i)_{i=1}^{M_n}$ and moreover if $|F|\leq M_n$ and $(a_j)_{j \in F}$
are scalars then

\begin{equation}\label{eq:smeq}
 \frac{1}{2}\|\sum_{j \in F} a_j \tilde{x}_j^{M_n} \|
 \leq   \|\sum_{j \in F} a_j x_j^{M_n} \|
 \leq 2 \|\sum_{j \in F} a_j \tilde{x}_j^{M_n} \|.
\end{equation}

\noindent By Remark~\ref{rmk:regular2} by passing to subsequences and relabeling, assume that
$(x_j^{M_n})_{(n,j) \in I}$ forms a regular array with basis constant at most equal to $2$.
Apply Theorem~\ref{Th:ASU}
to $(x_j^{M_n})_{(n,j) \in I}$ to get a subarray which satisfies the conclusion of
Theorem~\ref{Th:ASU}.  By relabeling call $(x_j^{M_n})_{(n,j) \in I}$ the resulting subarray.
Define

$$ y_j = \sum_{n=1}^j \delta_n x_{\ell_j}^{M_n}$$

\noindent where $(\ell_j)$ is an increasing sequence of positive integers which guarantees that
$|x_i^*y_i| < \varepsilon$.  Note that $(y_j)$ is weakly null since $(x_j^{M_n})_j$ is weakly
null for all $n$ and $\delta_n$ is summable.  Since $(x_j^{M_n})_{(n,j) \in I}$ is regular, $(y_j)$ is a basic sequence with
\begin{equation} \label{eq:5star}
\|y_j\| \geq \frac{\delta_1}{2} \inf_{n,m} \|x_m^n\|
\end{equation}
\noindent (since the basis constant of $(x_j^{M_n})_{(n,j) \in I}$ is at most equal to 2 by 
Remark~\ref{rmk:regular2}).
Since $(\delta_n)$ is summable, $(y_j)$ is also bounded.  Fix $n \in \N$ and let $n \leq F \subseteq \N$, with $|F|\leq M_n$.  Then

\begin{equation*}
\begin{split}
\| \sum a_j y_j \| \geq & \frac{1}{3}\| \sum_{j\in F; j\geq n} \delta_n a_j x_{\ell_j}^{M_n} \| \text{ (by Theorem~\ref{Th:ASU})}\\
 \geq & \frac{1}{6}\| \sum_{j\in F; j\geq n} \delta_n a_j \tilde x_j^{M_n} \| \text{ (by (\ref{eq:smeq}))}\\
 = & \frac{1}{6}\| \sum_{j\in F; j\geq n} \delta_n a_j \tilde x_{k_j}^{M_n} \| \\
   & \text{ where the map }     F \ni j \mapsto k_j \in \{1,2,\ldots, |F|\} \text{ is a 1-1 increasing function} \\
  \geq & \frac{1}{6C} \delta_n \| \sum_{j\in F; j\geq n} a_jz_j \| \text{ (by (\ref{eq:condB}))}
\end{split}
\end{equation*}

\noindent Thus
$$\|\sum a_jy_j \| \geq  \sup_n
        \frac{1}{6C} \sup_{n \leq F \subseteq \N; |F|\leq M_n} \delta_n \|\sum_{i \in F} a_i z_i\| .$$
\end{proof}

Theorem~\ref{Thm:3.1}, just as Theorem~\ref{Thm:main},  gives sufficient conditions
on a Banach space $X$, in order that the ``multiple of the inclusion plus compact'' 
problem has an affirmative
solution on $X$. If in Theorem~\ref{Thm:3.1} one considers the special case where 
$(\underline{x}_i)_i=(x_i)_i$,
then the assumptions of Theorem~\ref{Thm:3.1} are similar to the assumptions of
Theorem~\ref{Thm:main}. The difference in that case is that in Theorem~\ref{Thm:3.1} 
(but not in Theorem~\ref{Thm:main}),
we assume that the sequence $(x_i)$ satisfies Property P2. Instead, in Theorem~\ref{Thm:main} we 
assume
that the sequence $(z_i)$ has a spreading model which has Property P2.
In Proposition~\ref{P:P2} we show that if the basic sequence $(z_i)$ has a spreading 
model which has Property P2 then
by replacing $(z_i)$ by a new sequence and relabeling, we can assume that $(z_i)$ has a spreading 
model which is unconditional and has Property P1.
Then, in Proposition~\ref{P:P1usm} we show that if the basic sequence $(z_i)$ has a spreading model
which is unconditional and has Property P1 then the sequence
$(x_i)$  can be ``replaced'' by a sequence $(x_i^+)$ which (may not be contained in the Banach space 
$X$ and) satisfies Property P2. 
Thus the proof of Theorem~\ref{Thm:main} will follow from Theorem~\ref{Thm:3.1} and
Propositions~\ref{P:P2} and \ref{P:P1usm}.

\begin{Prop} \label{P:P2}
Let $(z_i)$ be a seminormalized basic sequence which has a spreading model which has Property P2.
Then there exists a subsequence $(z_{k_i})$ of $(z_i)$ such that $(Z_i)$ has a spreading model
which is unconditional and has Property P1, where either $(Z_i)_i:=(z_{k_i})_i$ or  
$(Z_i)_i:=(z_{k_{2i}}-z_{k_{2i-1}})_i$.
\end{Prop}

\begin{proof}
By Rosenthal's $\ell_1$ Theorem \cite{R2} there exists a subsequence $(z_{k_i})$ of $(z_i)$ such that
either $(z_{k_i})$ is equivalent to the standard basis of $\ell_1$, or $(z_{k_i})$ is weak Cauchy.
In the first case by passing to a further subsequence and relabeling assume that 
$(z_{k_i})$ has a spreading model and set $(Z_i)_i:=(z_{k_i})_i$. Then obviously
$(Z_i)$ has a spreading model which is unconditional and has Property P1. 
If $(z_{k_i})$ is weak Cauchy,
set $(Z_i)_i=(z_{k_{2i}}-z_{k_{2i-1}})_i$. Then $(Z_i)$ is weakly null, hence by \cite{BS1}, \cite{BS2}
we can pass to a subsequence of $(Z_i)$ and relabel in order to assume that $(Z_i)$ has a suppression
$1$-unconditional spreading model. It is obvious to see that the Property P2 passes from the 
spreading model of $(z_i)$ to the spreading model of $(Z_i)$.
Thus $(Z_i)$ has a spreading model which is unconditional and has Property P1 
(see Proposition~\ref{P:P1P2}(b)).
\end{proof}

\begin{Prop} \label{P:P1usm}
Let $(x_i)$ and $(z_i)$ be two seminormalized basic sequences (not necessarily in the same Banach 
space)
such that $(x_i)$ satisfies condition (\ref{eq:condA})  of Theorem~\ref{Thm:main}
and $(z_i)$ has a spreading model which is unconditional and has Property P1. Then there exists a 
seminormalized basic sequence
$(x_i^+)$ which has Property P2, dominates $(x_i)$ and has spreading model $(\tilde{x}^+_i)$ which
is s.c. dominated by $(z_i)$.
\end{Prop}

\begin{proof}[Proof of Proposition~\ref{P:P1usm}]
Before defining $(x^+_i)$, define an auxiliary basic sequence $(z_i')$ as follows. Define the norm
on the span of $(z_i')$ as the completion of the following:
for $(a_i) \in c_{00}$ let

\begin{equation} \label{eq:zprime}
  \|\sum a_iz_i' \| := \sup_{n \in \N}  \{ \frac{1}{\sqrt{L_n}}\|\sum_{i\in A} a_i \tilde z_i\|:
  A \subseteq \N \text{ with }|A|\leq n \}
\end{equation}

\noindent where
\begin{equation} \label{eq:Ln}
L_n = \sup_{1\leq k \leq n} \|\sum_{i =1}^k \tilde z_i\|.
\end{equation}

\noindent {\bf Claim 1:} The sequence  $(z_i')$ is seminormalized, $1$-spreading, unconditional
(thus $1$-subsymmetric) and has Property P1.

Since $(\tilde z_i)$
is $1$-spreading, the Property P1 of $(\tilde z_i)$ is equivalent to the fact that 
$L_n \rightarrow \infty$.
Also it is easy to verify that the $1$-spreading and unconditionality properties pass from 
$(\tilde z_i)$ to $(z_i')$.
Since $(L_n)$ is increasing we have that $\|z_i'\| = \frac{\|\tilde z_i\|}{\sqrt{L_1}} = 
\sqrt{\|\tilde{z_i}\|}$
therefore $(z_i')$ is seminormalized.  Notice that
$$
L_n = \sup_{1\leq k \leq n} \|\sum_{i =1}^k \tilde z_i\| \leq C_1 \|\sum_{i=1}^n \tilde z_i\|
$$
\noindent where $C_1$ is the basis constant of $(\tilde z_i)$.  Thus
$$
\|\sum_{i=1}^n z_i' \| \geq \frac{1}{\sqrt{L_n}} \|\sum_{i=1}^n \tilde z_i \|
                \geq \frac{1}{C_1}\sqrt{L_n}  \rightarrow \infty.
                $$
\noindent Hence $(z_i')$ has Property P1 (since $(z_i')$ is $1$-spreading). This finishes the proof of
Claim~1.

\noindent{\bf Claim 2:} $(z_i)>>(z_i')$.

Let $\varepsilon > 0$ be given.  We will choose $\rho >0$ so that
$\Delta_{(z_i),(z_i')}(\rho) \leq \varepsilon.$  Since $L_n \rightarrow \infty$
we can find an $N \in \N$ such that $C_2/\sqrt{L_N} < \varepsilon$
where $C_2$ is the suppression unconditionality constant of the sequence $(\tilde z_i)$
(by \cite{BS1,BS2} we have $C_2=1$ if $(z_i)$ is weakly null).  Let
$\rho = \min_{n \leq N} \frac{\varepsilon\sqrt{L_n}}{\|\tilde z_1\|n}$.  Let $(a_i) \in c_{00}$ be
such that  $\|\sum a_iz_i \|= 1$ and $|a_i| \leq \rho$.  Also let
$n_0 \in \N$ and $A$ be a subset of $\N$ with $|A| \leq n_0$ and
$\|\sum a_iz_i' \|= \frac{1}{\sqrt{L_{n_0}}}\|\sum_{i \in A} a_i \tilde z_i \|$.
If $n_0 \leq N$ then $\frac{1}{\sqrt{L_{n_0}}}\|\sum_{i \in A} a_i \tilde z_i \|
\leq \frac{\|\tilde z_1\|}{\sqrt{L_{n_0}}}\rho n_0 \leq \varepsilon$ by the choice of $\rho$
(notice that $\|\tilde z_i\| = \|\tilde z_1\|$ for all $i$).
If $n_0 > N$ then $\frac{1}{\sqrt{L_{n_0}}}\|\sum_{i \in A} a_i \tilde z_i \|
\leq \frac{1}{\sqrt{L_{N}}}C_2 \|\sum a_i \tilde z_i\| = \frac{C_2}{\sqrt{L_{N}}}< \varepsilon$
where the first inequality follows by the unconditionality of $(\tilde z_i)$.
Thus $\Delta_{(z_i),(z_i')}(\rho)<\varepsilon$.  This finishes the proof of Claim~2.

Now we are ready to define the basic sequence $(x_i^+)$.  Define the norm on the span of the
basic sequence $(x_i^+)$ as the completion of the following:  for $(a_i) \in c_{00}$ let
\begin{equation} \label{eq:xprime}
\|\sum a_ix_i^+ \| = \max \{\|\sum a_ix_i \|, \|\sum a_i z_i'\|  \}
\end{equation}
\noindent where $(z_i')$ is defined by (\ref{eq:zprime}).

\noindent{\bf Claim 3:} The sequence  $(x_i^+)$ is seminormalized, dominates $(x_i)$, has 
Property P2
and has spreading model $(\tilde x_i^+)$ which satisfies $(z_i)>>(\tilde x_i^+)$.

Indeed by Claim~1 we have that $(z_i')$ is seminormalized thus $(x_i^+)$ is seminormalized.  
Obviously, by
(\ref{eq:xprime}), we have that $(x^+_i)$ dominates $(x_i)$. By Claim~1
we have that $(z_i')$ is unconditional and has Property P1 thus  it is easy to see that
$(z_i')$ has Property P2, (see Proposition~\ref{P:P1P2} (c)).  Since $(x_i^+)$
dominates $(z_i')$, we obtain that $(x_i^+)$ has Property P2.  Since $(x_i)$ has spreading model 
$(\tilde x_i)$
and $(z_i')$ is $1$-spreading, (\ref{eq:xprime}) implies that $(x_i^+)$ has spreading model 
$(\tilde x_i^+)$
which satisfies
\begin{equation} \label{eq:21}
\|\sum a_i \tilde x_i^+\| = \max \{ \|\sum a_i \tilde x_i\|, \|\sum a_i z_i'\| \}.
\end{equation}
\noindent  Finally, Claim~2, the fact that $(z_i) >> (\tilde x_i)$ and (\ref{eq:21}) imply that
$(z_i) >> (\tilde x_i^+)$. This finishes the proof of Claim~3 and of Proposition~\ref{P:P1usm}.
\end{proof}

Now we are ready to present the proof of Theorem~\ref{Thm:main}.

\begin{proof}[Proof of Theorem \ref{Thm:main}]
We know that $(z_i)$ has a spreading model which has Property P2. Then by Proposition~\ref{P:P2},
there exists a subsequence $(z_{k_i})$ of $(z_i)$ such that $(Z_i)$ has a spreading model which is
unconditional and has Property P1, where either $(Z_i)_i:=(z_{k_i})_i$ or  
$(Z_i)_i:=(z_{k_{2i}}-z_{k_{2i-1}})_i$. 

First assume that  $(Z_i)_i=(z_{k_i})_i$. Then set $(X_i)_i:=(x_i)_i$ and $(X^n_i)_i:=(x^{k_n}_i)_i$. 
We claim that 
(\ref{eq:condA}) and (\ref{eq:condB}) are satisfied with ``$(x_i)_i$'', ``$(x^n_i)_i$'' and ``$(z_i)_i$'' being 
replaced by $(X_i)_i$, $(X^n_i)_i$ and $(Z_i)$ respectively. Indeed, since $(\tilde{x}_i)<<(z_i)$, 
we have that $(\tilde{x}_{k_i})<<(Z_i)$. 
Since $(\tilde{x}_i)$
is isometrically equivalent to $(x_{k_i})$, we obtain that $(\tilde{x}_i) <<(Z_i)$. Thus (\ref{eq:condA}) is
satisfied for $(X_i)_i$ and $(Z_i)_i$.
Also notice that (\ref{eq:condB}) is satisfied with ``$(x^n_i)$'' and ``$(z_i)$'' being replaced by 
$(X^n_i)$ and $(Z_i)$
respectively. Indeed, since $(z_i)_{i=1}^n$ is $C$-dominated by $(\tilde{x}^n_i)$, we have that
$(Z_i)_{i=1}^n$ is $C$-dominated by $(\tilde{x}^{k_n}_{k_i})_{i=1}^n$ which is isometrically equivalent
to $(\tilde{x}^{k_n}_i)_{i=1}^n$. 

Second assume that $(Z_i)_i:=(z_{k_{2i}}-z_{k_{2i-1}})_i$. Then set $(X_i)_i:=(x_{2i}-x_{2i-1})_i$ and
$(X^n_i)_i:=(x^{k_{2n}}_{2i}-x^{k_{2n}}_{2i-1})_i$. We claim that (\ref{eq:condA}) and (\ref{eq:condB})
are satisfied for ``$(x_i)$'', ``$(x^n_i)$'' and ``$(z_i)$'' being replaced by $(X_i)$, $(X^n_i)$ and 
$(Z_i)$ respectively. Indeed notice that since $(\tilde{x}_i) <<(z_i)$, we have that
$(\tilde{x}_{k_{2i}} -\tilde{x}_{k_{2i-1}}) <<(Z_i)$. Since $(\tilde{x}_{k_{2i}} -\tilde{x}_{k_{2i-1}})$ is
isometrically equivalent to $(\tilde{x}_{2i} -\tilde{x}_{2i-1})$ and $(X_i)$ has spreading model
isometrically equivalent to $(\tilde{x}_{2i} -\tilde{x}_{2i-1})$, we have that (\ref{eq:condA}) is satisfied
for $(X_i)$ and $(Z_i)$. Also, since $(z_i)_{i=1}^n$ is $C$-dominated by $(\tilde{x}^n_i)_{i=1}^n$,
we have that $(Z_i)_{i=1}^n$ is $C$-dominated by 
$(\tilde{x}^{k_{2n}}_{k_{2i}}-\tilde{x}^{k_{2n}}_{k_{2i-1}})_{i=1}^n$,
which is isometrically equivalent to $(\tilde{x}^{k_{2n}}_{2i}-\tilde{x}^{k_{2n}}_{2i-1})_{i=1}^n$. Finally
notice that the spreading model of $(X^n_i)_i$ is isometrically equivalent to
$(\tilde{x}^{k_{2n}}_{2i}-\tilde{x}^{k_{2n}}_{2i-1})_i$. Therefore condition (\ref{eq:condB}) is satisfied
for $(X^n_i)_i$ and $(Z_i)_i$.

In either of the above two cases (i.e. either $(Z_i)_i=(z_{k_i})_i$ or $(Z_i)_i:=(z_{k_{2i}}-z_{k_{2i-1}})_i$),
since $(\tilde{X}_i) <<(Z_i)$ and $(Z_i)$  has a spreading model which is unconditional and has 
Property P1, by Proposition~\ref{P:P1usm},
there exists a seminormalized basic sequence
$(x_i^+)$ which has Property P2, dominates $(X_i)$ and has spreading model $(\tilde{x}^+_i)$ which
is s.c. dominated by $(Z_i)$. Apply Theorem~\ref{Thm:3.1} for ``$(\underline{x}_i)_i$'' 
being equal to $(X_i)$, ``$(x^n_i)$'' being equal to $(X^n_i)$, ``$(z_i)$'' being equal to ``$(Z_i)$'' 
and  ``$(x_i)_i$'' being equal to $(x^+_i)$ to finish the proof.
\end{proof}

\begin{proof} [Proof of Theorem \ref{thm:AOST_p}]
If $p$ belongs to the Krivine set of $(\tilde{x_i})$ then for all $n \in \N$ there
exists $(x_i^n)_{i \in \N}$ a block sequence of $(x_i)$ of identically
distributed blocks such that any $n$ terms of $(x_i^n)_{i \in \N}$ are
$2$-equivalent to the unit vector basis of $\ell_p^n$.
Then apply Theorem~\ref{Thm:main} for $(z_i)$ being the unit vector basis of $\ell_p$.
\end{proof}

Next we examine the relation between Properties P1 and P2 and how these properties pass to spreading
models.

\begin{Prop} \label{P:P1P2}
\begin{itemize}
\item[(a)] Assume that a seminormalized basic sequence $(z_i)$ has a spreading model $(\tilde{z}_i)$.
If $(z_i)$ has Property P1 then $(\tilde{z}_i)$ has Property P1. If $(z_i)$ has Property P2 then
$(\tilde{z}_i)$ has Property P2.
\item[(b)] If a seminormalized basic sequence has Property P2 then it has Property P1.
\item[(c)] If an unconditional seminormalized basic sequence has Property P1 then it has Property P2.
\item[(d)] Let $(z_i)$ be a Schreier unconditional seminormalized basic sequence in a Banach space.
        If $(z_i)$ has Property P1 then there exists some subsequence $(z_{n_i})$ of $(z_{i})$ which has Property P2.
\end{itemize}
\end{Prop}

\begin{proof} Parts (a), (b) and (c) are obvious.
For part (d) assume (towards contradiction) that no subsequence of $(z_i)$ has Property P2.  So
for every subsequence $(m_i)$ of $\N$ there exists $\rho>0$ such that for all $M \in \N$ there exists
$(a_i) \in c_{00}$ with $\|\sum a_iz_{m_i}\| = 1$ and $|\{i:|a_i|>\rho \}|>M$.  Thus
$\rho$ depends on the sequence $(m_i)$.  Let $\rho((m_i))$ be the supremum of such $\rho$'s.
Notice that if $(m_i)$, $(n_i)$ are subsequences of $\N$ with
$(n_i) \subseteq (m_i)$ then $\rho((n_i)) \leq \rho((m_i))$.  We claim there is a subsequence
$(m_i)$ of $\N$ such that $\inf \{ \rho((m_{n_i}) : (n_i) \text{ subsequence of } \N \}>0$.
To show this claim we assume again by contradiction that for all subsequences $(m_i)$ of $\N$ we have
$\inf \{ \rho((m_{n_i})) : (n_i) \text{ subsequence of } \N \} = 0$.

Let $(n_i^1)_{i=1}^\infty \subseteq \N$ be such that $\rho((n_i^1)) < 1.$
Let $(n_i^2)_{i=1}^\infty \subseteq (n_i^1)_{i=1}^\infty$ be such that
$\rho((n_i^2)) < \frac{\rho((n_i^1))}{2C_s}$, where $C_s$ is the constant of Schreier unconditionality of
$(z_i)$.  Given $(n_i^k)_{i=1}^\infty$ define $(n_i^{k+1})_{i=1}^\infty \subseteq (n_i^k)_{i=1}^\infty$
to be a subsequence such that $\rho((n_i^{k+1})) < \frac{\rho((n_i^k))}{2C_s}$.  Define $n_i := n_i^i$.
By our assumption $\rho((n_i)) > 0$.  Let $k \in \N$ be such that $\rho((n_i)) >  \rho((n_i^k)) > 0$.
Let $N \in \N$ be arbitrary such that $N > k$.   By the definition of $\rho((n_i))$ there exist a sequence of scalars
$(a_{n_i})$ such that $\|\sum a_{n_i}z_{n_i}\| = 1$ and $|\{n_i: |a_{n_i}| > \frac{\rho((n_i))}{2} \}| > 2(N+1).$
Define $A := \{n_i: |a_{n_i}|>\frac{\rho((n_i))}{2} \} = \{\ell_1, \ell_2, \ldots ,\ell_{|A|} \}$ where
$\ell_1<\ell_2<\cdots$.  Define $m = \lfloor \frac{|A| + 1}{2}\rfloor$ and
$B =  \{\ell_{m+1}, \ell_{m+2}, \ldots ,\ell_{|A|} \}$.  Notice that $B$ is a Schreier subset of the sequence
$(n_i^{k+1})$ and $|B| \geq N $.  Thus we can project to the set $B$:
$$\|\sum_{n_i \in B} a_{n_i}z_{n_i}\| \leq C_s \|\sum a_{n_i}z_{n_i}\| = C_s .$$
\noindent Thus $\|\sum_{n_i \in B} \frac{a_{n_i}}{C_s}z_{n_i}\| \leq 1$.   Note that
$B \subseteq \{n_i: \frac{|a_{n_i}|}{C_s}>\frac{\rho((n_i))}{2C_s} \}$.  Since $B \subset (n_i^{k+1})$ we have that
$\rho((n_i^{k+1})) \geq \frac{\rho((n_i))}{2C_s}$.  But $\rho((n_i)) > \rho((n_i^k)) > 2C_s \rho((n_i^{k+1}))$, a
contradiction. Thus there exists some subsequence $(n_i)$ where
$\inf \{ \rho((n_{m_i})) : (m_i) \text{ subsequence of } \N \} > 0$.

For ease of notation call $(z_i)$ the subsequence $(z_{n_i})$.
Let $\rho > 0$ such that
\begin{equation} \label{eq:pedro}
\begin{split}
& \text{ for all sequences } (s_i) \text{ of positive integers and } m \in \N \text{ there exists } \\
& (a_i) \in c_{00}\text{ with }\|\sum a_{s_i} z_{s_i}\| = 1 \text{ and } |\{i : |a_{s_i}|>\rho \}|>m.
\end{split}
\end{equation}

Fix $m \in \N.$ Let
\begin{equation*}
\begin{split}
\sigma_m = \{(\ell_i)_{i=1}^m &\subseteq \N:\text{ there exists }(a_i) \in c_{00}
        \text{ such that }\|\sum a_iz_i\| = 1 \\
 & \text{ and } (\ell_i)_{i=1}^m \subseteq \{i:|a_i|>\rho \} \}.
\end{split}
\end{equation*}
\noindent By Ramsey's theory there exists some
$(m_i)_{i=1}^\infty \subset \N$ such that either $[(m_i)_{i=1}^\infty]^m \subset \sigma_m$ or
$[(m_i)_{i=1}^\infty]^m \cap \sigma_m = \emptyset$.  But the second case is not possible by (\ref{eq:pedro}).

Thus for each $m \in \N$  there exists some subsequence $(n_i^m)$ of $\N$ such that for each
$F \subseteq (n_i^m)$ where $|F| = m$, there exists some $(a_i) \in c_{00}$ such that
$\|\sum a_i z_i\| = 1 $ and $F \subseteq \{i:|a_i|>\rho\}.$  We can choose each of these
subsequences such that $(n_i^1) \supseteq (n_i^2)\supseteq (n_i^3) \supseteq \cdots$ and then
define $n_i = n_i^i$.  Thus for all $A \subseteq \{m,m+1,\ldots\} $ such that $|A|=m$
(i.e. $A$ is Schreier) there exists some $(a_i^m) \in c_{00}$ such that
$\|\sum a_i^mz_i\| = 1$ and $(n_i)_{i \in A} \subseteq \{i: |a_i^m|>\rho\}$.  So we define for each
$m \in \N$ such $(a_i^m)$ where $\|\sum a_i^mz_i\| = 1$ and $(n_i)_{i =m}^{2m-1} \subseteq \{i: |a_i^m|>\rho\}$.
Notice each $a_i^m \in \mathcal{A} := \{z \in \mathbb{C}: \rho < |z| < \frac{2C}{\inf_i \| z_i \|} \}$
where $C$ is the basis constant of $(z_i)$. It is an elementary exercise to see the following.

\begin{equation} \label{eq:elem}
\begin{split}
& \text{For all }M \in \N \text{ and }\varepsilon >0 \text{ there exists }k \in\N \text{ such that if }
                (a_i)_{i=1}^k \subseteq \mathcal{A}\\
&  \text{ then there exists some subset } A' \subseteq \{1,2,\dots,k\}
        \text{ such that }\\
& |A'| \geq M \text{ and for all } i,j \in A' \text{ we have } |a_i-a_j|<\varepsilon.
\end{split}
\end{equation}

Since $(z_i)$ has Property P1,
$$\liminf_{n \rightarrow \infty} \inf_{A \subseteq \N; |A|=n}\|\sum_{i \in A} z_i\| = \infty.$$
\noindent So there exists $N \in \N$ such that for all sets $A \subseteq \N$ with $|A|=N$ we have
\begin{equation} \label{eq:big}
\|\sum_{i\in A} z_i\| > \frac{C_s}{\rho}.
\end{equation}

By (\ref{eq:elem}) for $\varepsilon = \frac{C_s}{2N}$ and $M = N$ to get $k \in \N$.  Then note that
$(a_i^k)_{i=k}^{2k-1} \subseteq \mathcal{A}$ has $k$ many terms (some terms may be equal).
Therefore there is an $N$ element set
$A' \subseteq \{k, k+1, \ldots ,2k-1 \}$ such that for all $i, j \in A'$ we have
$|a_{i}^k-a_{j}^k| < \varepsilon$ .  Notice that $A'$ is
a Schreier set.  Thus by Schreier unconditionality we have

\begin{equation*}
\begin{split}
C_s = & C_s\|\sum a_i^k z_i\| \geq \|\sum_{i \in A'} a_i^k z_i\| \geq
                \|\sum_{i \in A'} a_1^k z_i\|  - \|\sum_{i \in A'} (a_i^k-a_1^k) z_i\| \\
        &\geq \|\sum_{i \in A'} a_1^k z_i\|  - N \max |a_i^k-a_1^k| \geq \|\sum_{i \in A'} a_1^k z_i\|  - \frac{C_s}{2}.
\end{split}
\end{equation*}

\noindent Therefore $\|\sum_{i \in A'} z_i\| \leq \frac{C_s}{2a_1^k} \leq \frac{C_s}{2\rho}$ which is a
contradiction to (\ref{eq:big}).
\end{proof}

Note that the summing basis has Property P1 but no subsequence of it has Property P2. Thus the 
assumption
that $(z_i)$ is Schreier unconditional in Proposition~\ref{P:P1P2}(d) is needed.

Notice that if $(x_i)$, $(x^n_i)_i$ (for $n \in \N$) and $(z_i)$ are seminormalized basic sequences
satisfying conditions (\ref{eq:condA}) and (\ref{eq:condB}) of Theorem~\ref{Thm:main}, and 
$(z_{k_i})$ is any subsequence of $(z_i)$ then $(x_i)$, $(x^{k_n}_i)_i$ (for $n \in \N$) and 
$(z_{k_i})$ also satisfy conditions (\ref{eq:condA}) and (\ref{eq:condB}). Indeed, the sequence 
$(\tilde{x}_i)$ is isometrically equivalent to the sequence $(\tilde{x}_{k_i})$ which is s.c. dominated 
by $(z_{k_i})$. Also $(z_{k_i})$ is $C$-dominated by $(\tilde{x}^{k_n}_{k_i})_{i=1}^n$ which is 
isometrically equivalent to $(\tilde{x}^{k_n}_i)_{i=1}^n$. 
This observation, Theorem~\ref{Thm:main} and Proposition~\ref{P:P1P2}(a) and (c), immediately 
give the following result.

\begin{Cor}
Let $X$ be a Banach space containing seminormalized basic sequences $(x_i)_i$ and $(x^n_i)_i$
for all $n \in \N$, such that $0 < \inf_{n,i} \| x^n_i\| \leq \sup_{n,i} \| x^n_i \| < \infty$.
Let  $(z_i)_i$ be a basic sequence not necessarily
in $X$. Assume that $(x_i)$ satisfies condition (\ref{eq:condA}) 
and for all $n \in \N$ the sequence $(x^n_i)_i$ satisfies condition (\ref{eq:condB}) of 
Theorem~\ref{Thm:main}.
Assume that the sequence $(z_i)$ satisfies at least one of the following three conditions:
\begin{itemize} 
\item a subsequence of $(z_i)$ has  a spreading model which has Property P2, or
\item a subsequence of $(z_i)$ has Property P1 and has a spreading model which is unconditional, or
\item a subsequence of $(z_i)$ has a spreading model which is unconditional and has Property P1, or
\item a subsequence of $(z_i)$ has Property P2. 
\end{itemize}
Then there exists a subspace $Y$ of $X$ which has a basis and an operator
$T \in \mathcal{L}(Y,X)$ which is not a compact perturbation of a multiple of the inclusion map.
\end{Cor}

%------------------------------------------------------------------------------------------------
\section{An Application of Theorem~\ref{Thm:main}} \label{sec:dew}
%------------------------------------------------------------------------------------------------

Next we give an application of Theorem~\ref{thm:AOST_p} where previously known results do not seem
to be applicable (at least with the same ease).  As mentioned before the problem of finding a subspace
$Y$ of a Banach space $X$ and an operator $T \in \mathcal{L}(Y,X)$ which is not a compact perturbation
of the inclusion operator is non trivial when $X$ is saturated with HI Banach spaces.  The
HI space to which Theorem~\ref{thm:AOST_p} will be applied was constructed by
N.~Dew \cite{D}, and here will be denoted by $D$.  The construction of the space $D$ is based on the
$2$-convexification of the Schlumprecht space $S$ \cite{S1} in a similar manner that the space of T.W.~Gowers
and B.~Maurey \cite{GM} is based on $S$.  We recall the necessary definitions.

Let $X$ be a Banach space with a basis $(e_i)$.  For any interval $E$ in $\N$
and a vector $x = \sum x_je_j \in X$ define $Ex = \sum_{j \in E} x_je_j \in X$.
There is a unique norm
$\|\cdot\|_S$ on $c_{00}$ which satisfies:

$$\|x\|_S = \sup \left \{\frac{1}{f(\ell)}\sum_{i=1}^\ell \|E_ix\|_S:
        E_1<\cdots<E_\ell \right \} \vee  \|x\|_{\ell_\infty}$$

\noindent  where $f(\ell) = \log_2 (\ell +1)$.   The completion of $c_{00}$ under this norm is the
Banach space $S$.  Let $S_2$ be its 2-convexification.   Recall if $X$ is a Banach space
having an unconditional basis $(e_n)$, then we can define the $2$-convexification
$X_2$ of $X$ by the norm

$$\|\sum a_n \sqrt{e_n} \|_{X_2} := (\|\sum a_n^2e_n\|_{X})^{1/2}$$
where $(\sqrt{e_n})$ denotes the basis of $X_2$, (we will talk more later about the ``square root'' map 
from $X$ to $X_2$). 
We will show that the spreading model of the
unit vector basis of $D$ is the unit vector basis of $S_2$.  Before we do this we must see
the definition of $D$.

In order to define the Banach space $D$, a lacunary set $J \subseteq \N$ is used which has the 
property that
if $n,m \in J$ and $n<m$ then $8n^4 \leq \log \log \log m$, and $f(\min J)
\geq 45^4$.  Write $J$ in increasing order as $\{ j_1, j_2, \ldots \}$.
Now let $K \subset J$ be the set $\{ j_1, j_3, j_5 ,
\cdots \}$ and $L \subset J$ be the set $\{ j_2 , j_4 , j_6 , \ldots
\}$. Let ${\bf Q}$ be the set of scalar sequences
with finite support and rational coordinates  whose absolute value is at most one.
Let $\sigma$ be an injective function from ${\bf Q}$ to $L$ such that if
$z_1, \ldots , z_i$ is such a sequence, then $(1/400)f(\sigma(z_1,
\ldots , z_i)^{1/40})^{\frac{1}{2}} \geq \# \supp (\sum_{j=1}^i z_j)$.
Then, recursively, we define a set
of functionals of the unit ball of the dual space $D^*$ as follows: Let
$$
D^*_0= \{ \lambda e_n^* : n \in \N, \ | \lambda | \leq 1 \} .
$$
Assume that $D_k^*$ has been defined. Define the norm $\| \cdot \|_k$ on $c_{00}$
by

\begin{equation} \label{eq:normk}
\| x \|_k = \sup \{|x^*(x)|: x^* \in D_k^* \}
\end{equation}

\noindent and let $\| \cdot \|_k^*$ denote its dual norm.  Then $D_{k+1}^*$ is the set of
all functionals of the form $F \ z^*$ where $F \subseteq \N$ is an interval
and $z^*$ has one of the following three forms:
\begin{equation} \label{eq:f1}
z^*= \sum_{i=1}^\ell \alpha_i z_i^*
\end{equation}
\hskip .4in where $\sum_{i=1}^\ell | \alpha_i | \leq 1$ and $z_i^* \in
D_k^*$  for $i=1, \ldots , \ell$.
\begin{equation}\label{eq:f2}
z^*= \frac{1}{\sqrt{f(\ell)}} \sum_{i=1}^\ell \alpha_i z_i^*
\end{equation}
\hskip .4in where $\sum_{i=1}^\ell  \alpha_i^2  \leq 1$,
$z_i^* \in D_k^*$ for $i=1, \ldots , \ell$, and
$z_1^*  < \cdots <z^*_\ell$.

\begin{equation} \label{eq:f3}
z^*= \frac{1}{\sqrt[4]{f(\ell)}} \sum_{i=1}^\ell z_i^*
\text{ where }
z_i^* = \frac{1}{\sqrt{f(m_i)}} \sum_{j=1}^{m_i}
\alpha_{i,j} \frac{Ez_{i,j}^*}{\|Ez_{i,j}^*\|_{k}^*}
\end{equation}
\hskip .4in for certain $(\alpha_{i,j})$ where $\sum_{i,j}
\alpha_{i,j}^2 \leq 1$ ($\alpha_{i,j}$'s are explicitly
chosen in \cite{D}, but the exact values are not needed for our
purposes) where $z^*_{i,j} \in D^*_k$ for $1 \leq i \leq \ell$ and $1 \leq j
\leq m_i$, $z^*_{1,1} < \cdots < z^*_{1,m_1}< z^*_{2,1} < \cdots <
z^*_{\ell , m_\ell}$, $m_1= j_{2 \ell}$, $\beta z_i^*$ has rational
coordinates for some $\beta >0$ (whose exact value is not needed for our purpose),
$m_{i+1}= \sigma (\beta z_1^*,\ldots , \beta z_i^*)$, for $i=1, \ldots , \ell -1$
and $E$ is an interval.

Finally, the norm of $D$ is defined by
$$
\| x \|_{D} = \sup \{ z^* (x) : z^* \in \cup_{k=0}^\infty D^*_k \}.
$$

\begin{Prop} \label{prop:spreading}
The spreading model of the unit vector basis of $D$ is the unit
vector basis of $S_2$.
\end{Prop}

\begin{proof}
From the definitions of the two spaces it is easy to see that for $(a_i) \in c_{00}$,
$\| \sum a_i e_i \|_{S_2} \leq \| \sum a_i e_i \|_D$ ($(e_n)$ will denote the bases of both
spaces $S_2$ and $D$ but there will be no confusion about which space we consider at each
moment).  Thus to show Proposition~\ref{prop:spreading} we need only show
for any given $\varepsilon > 0$, and finitely many scalars $(a_i)_{i=1}^N$ there exists
$n_0$ such that for any $n_1, n_2, \ldots , n_N \in  \N$ with $n_0<n_1< n_2 < \cdots <n_N$, we have

$$\| x \|_D \leq \| y  \|_{S_2} + \varepsilon.$$

\noindent where $x = \sum_{i=1}^N a_ie_{n_i} \in D$ and $y = \sum_{i=1}^N a_ie_i \in S_2$.
This will follow immediately once we show by induction on $n$ that for any $n \in \N$,
$\varepsilon > 0$ and scalars $(a_i)_{i=1}^N$ we have

\begin{equation} \label{eq:induction}
\| x \|_n \leq \| y  \|_{S_2} + \varepsilon
\end{equation}

\noindent where $\| \cdot \|_n$ is defined in (\ref{eq:normk}).  For ``$n=0$''
we have $\|x\|_0 = \max_{1 \leq i \leq N} |a_i| \leq \|y\|_{S_2}$.  Now for the inductive step
assume that (\ref{eq:induction}) is valid for $n$.  Let $\varepsilon > 0$ and
scalars $(a_i)_{i=1}^N$.

First note that by the induction hypothesis there exists
$n_0^1 \in \N$ such that for all $n_0^1 <n_1<n_2< \cdots <n_N$ we have
\begin{equation} \label{eq:case1}
\| x \|_n \leq \| y  \|_{S_2} + \varepsilon .
\end{equation}

Secondly, by the inductive hypothesis there exists $n_0^2 \in \N$ such that for all 
$1 \leq i_0 \leq j_0 \leq N$
and $n_0^2 <n_1<n_2< \cdots <n_N$ we have
\begin{equation} \label{eq:case2}
\| \sum_{i=i_0}^{j_0} a_ie_{n_i} \|_n \leq \| \sum_{i=i_0}^{j_0} a_ie_i  \|_{S_2} + 
\frac{\varepsilon}{\sqrt{N}}.
\end{equation}

And finally there exists $j_0 \in L$ such that $\frac{1}{\sqrt{f(j_0)}}N \| y \|_{S_2} \leq \varepsilon $.  
Let $G = \sigma^{-1}(\{1,2,\ldots,j_0-1\})$, So $G$ is a finite subset of finite sequences of vectors 
with rational
coefficients.  Let $n_0^3$ be the maximum of the support of any vector in any sequence in $G$.  
And then set
$n_0 = \max \{n_0^1, n_0^2, n_0^3\}$.  Recall the norming vectors
$z^* \in D_{n+1}^*$ can be one of three different types.  Each type of functional
will present us with a different case.

\noindent {\bf CASE 1:}  Let $z^*$ be given by (\ref{eq:f1}).
Then

\begin{equation*}
\begin{split}
|z^*(x)| & \leq \sum_{i=1}^\ell |\alpha_i| |z_i^*(x)| \\
        & \leq \sum_{i=1}^\ell |\alpha_i| (\|x\|_{S_2} + \varepsilon) \text{ (by the (\ref{eq:case1}))} \\
        & = \|x\|_{S_2} + \varepsilon.
\end{split}
\end{equation*}

\noindent {\bf CASE 2:}  Let $z^*$ be given by (\ref{eq:f2}).
Thus for $n_0 <n_1<n_2< \cdots <n_N$ we have

\begin{equation*}
\begin{split}
|z^*(x)| &\leq \frac{1}{\sqrt{f(\ell )}} \sum_{j=1}^\ell |\alpha_j| |z_j^*(x)|  \\
 &\leq \frac{1}{\sqrt{f(\ell )}} \sum_{j=1}^\ell |\alpha_j| |z_j^*(E_jx)|
\end{split}
\end{equation*}
\noindent where $E_j$ is the smallest interval containing the support of $z_j^*$ intersected with the
support of $x$.  Continuing the above calculation we have

\begin{equation*}
\begin{split}
 |z^*(x)| &\leq \frac{1}{\sqrt{f(\ell )}} \sum_{j=1}^\ell |\alpha_j| \|E_jx\|_n \\
 &\leq \frac{1}{\sqrt{f(\ell )}} \sum_{j=1}^\ell |\alpha_j| (\|E_jx\|_{s_2} + \frac{\varepsilon}{\sqrt{N}})
        \text{ (by (\ref{eq:case2}))}\\
 &\leq \frac{1}{\sqrt{f(\ell )}} (\sum_{j=1}^\ell |\alpha_j|^2)^{1/2}
        (\sum_{j=1}^\ell \|E_jx\|_{s_2}^2)^{1/2} + \varepsilon\\
 &       \text{ (by Cauchy-Schwarz' inequality since there are at most $N$ many nonempty $E_j$'s )}\\
 &\leq \|x\|_{s_2} + \varepsilon.
\end{split}
\end{equation*}

\noindent {\bf CASE 3:}  Let $z^*$ be given by (\ref{eq:f3}).
We can of course assume $z^*(x) \not = 0$ thus let
$i_0$ be the smallest natural number such that

$$ \cup_{i=1}^{i_0} \supp(z_i^*) \cap \{n_0+1, n_0 +2, \ldots\} \not = \emptyset.$$

\noindent Then by the definition of $n_0^3$ we have that
$\sigma(\beta z_1^*,\beta z_2^*,\ldots,\beta z_i^*)>j_0$ for
$ i_0 + 1 \leq i \leq \ell $.  Thus

\begin{equation} \label{E:36}
\begin{split}
 & |z^*(x)| \\
 &\leq \frac{1}{\sqrt[4]{f(\ell )}} \bigg | \frac{1}{\sqrt{f(m_{i_0})}}\sum_{j=1}^{m_{i_0}} \alpha_{i_0,j}z_{i_0,j}^*(x) \bigg |
        + \frac{1}{\sqrt[4]{f(\ell)}}\frac{1}{\sqrt{f(m_{i_0+1})}}\sum_{j=1}^{m_{i_0+1}} |\alpha_{i_0+1,j}||z_{i_0+1,j}^*(x)|
        + \cdots \\
 &\leq \frac{1}{\sqrt[4]{f(\ell)}} \| y \|_{S_2} + \varepsilon
        + \frac{1}{\sqrt[4]{f(\ell)}}\frac{1}{\sqrt{f(m_{i_0+1})}}\sum_{ \{ 1 \leq j \leq m_{i_0+1}:z^*_{i_0+1,j}x \not = 0 \}} |\alpha_{i_0+1,j}||z_{i_0+1,j}^*(x)|
        + \cdots \text{ (by CASE 2)}\\
 &\leq \frac{1}{\sqrt[4]{f(\ell)}} \| y \|_{S_2} + \varepsilon
        + \frac{1}{\sqrt[4]{f(\ell)}}\frac{1}{\sqrt{f(j_0)}} N (\| y \|_{S_2} + \varepsilon),
\end{split}
\end{equation}
\noindent where the last inequality is valid since $|\alpha_{i,j}| \leq 1$, 
$|z^*_{i,j}(x)| \leq \| z_{i,j}^* \|^*_n \| x \|_n \leq \| x \|_n \leq \| y \|_{S_2} + \varepsilon$,
$m_i \geq j_0$ for all $i_0+1 \leq i \leq \ell$, (by the choice of $n_0$), 
 and there are at most $N$ many indices $(i,j)$ for which
$z_{i,j}^*(x) \not = 0$. The last expression in equation (\ref{E:36}) is at most equal to
 $\| y \|_{S_2} + 3\varepsilon $ which finishes the proof.
 \end{proof}

For the following remark we will need a bit more notation.  If $X$ is a Banach space
having an unconditional basis $(e_n)$ and $X_2$ is the $2$-convexification of $X$ then
for $x = \sum a_ne_n \in X$ then there is a canonical image of $x$ in $X_2$
which we define as  $\sqrt{x} = \sum \text{sign}(a_n) \sqrt{|a_n|}\sqrt{e_n} \in X_2$
(where $(\sqrt{e_n})$ denotes the basis of $X_2$).
\begin{Rmk} \label{rmk:krivine2}
Let $1 \leq p < \infty.$
Let $X$ be a Banach space with an unconditional basis $(e_i)_i$ and let $(\sqrt{e_i})_i$
be the basis of $X_2$. Then $p$ is in the Krivine set of $(e_i)$
if and only if $2p$ is in the Krivine set of $(\sqrt{e_i})$.
\end{Rmk}
\begin{proof}
Let $n \in \N$ and $\varepsilon > 0$ be given.  Since  $p$ is in the Krivine set of $X$, there exists
a block sequence  $(v_i)_{i=1}^n$ of $(e_i)$ such that for any
scalars $(a_i)_{i=1}^n$ we have that

$$ \frac{1}{1+\varepsilon}(\sum_{i=1}^n a_i^p)^{1/p} \leq \|\sum_{i=1}^n a_i v_i \| \leq
        (1+\varepsilon)(\sum_{i=1}^n a_i^p)^{\frac{1}{p}}$$

\noindent Let $w_i = \sqrt{v_i} \in X_2$, and scalars $(a_i)_{i=1}^n$ then

$$ \|\sum_{i=1}^n a_i w_i \|_{X_2} = \|\sum_{i=1}^n a_i^2 v_i \|_{X}^{\frac{1}{2}} \leq
        (1+\varepsilon)(\sum_{i=1}^n a_i^{2p})^{\frac{1}{2p}}$$

\noindent and
$$ \|\sum_{i=1}^n a_i w_i \|_{X_2} = \|\sum_{i=1}^n a_i^2 v_i \|_{X}^{\frac{1}{2}} \geq
        \frac{1}{1+\varepsilon}(\sum_{i=1}^n a_i^{2p})^{\frac{1}{2p}}.$$
The proof of the converse is similar.
\end{proof}

It is known \cite{S1} that the Krivine set of the unit vector basis of $S$
consists of the singleton $\{1\}$.  Thus by Remark~\ref{rmk:krivine2} we have:

\begin{Rmk} \label{rmk:star}
The Krivine set of the unit vector basis of $S_2$ consists of the singleton $\{2\}$.
\end{Rmk}

\begin{Prop} \label{prop:stillD}
Let $X$, $Y$ be Banach spaces with unconditional bases, $(x_n)$ be a basic sequence in $X$,
$(y_n)$ be a basic sequence in $Y$ such that $(x_n) >> (y_n)$.
Then $(\sqrt{x_n}) >> (\sqrt{y_n})$ where $(\sqrt{x_n})$ and $(\sqrt{y_n})$ are
the canonical images of $(x_n)$ and $(y_n)$ in $X_2$ and $Y_2$ respectively.
\end{Prop}
\begin{proof}
Note that

\begin{equation} \label{eq:33}
\begin{split}
\liminf_{n \rightarrow \infty} \inf_{A \subseteq \N; |A|= n} \| \sum_{i \in A} \sqrt{x_i} \|_{X_2} =
\liminf_{n \rightarrow \infty} \inf_{A \subseteq \N; |A|= n} \| \sum_{i \in A} x_i \|_{X}^{\frac{1}{2}}
= \infty.
\end{split}
\end{equation}

\noindent Also
\begin{equation} \label{eq:34}
\begin{split}
\Delta_{(\sqrt{x_n}),(\sqrt{y_n})}(\varepsilon) &= \sup\{\|\sum a_i\sqrt{y_i}\|_{Y_2}: |a_i|\leq \varepsilon
        \text{ and } \|\sum a_i \sqrt{x_i}\|_{X_2}=1 \} \\
  &= (\sup\{\|\sum a_i^2y_i\|_{Y}: |a_i|^2\leq \varepsilon^2 \text{ and }
        \|\sum a_i^2 x_i\|_{X}=1 \})^{\frac{1}{2}} \\
  &\leq (\Delta_{(x_n),(y_n)}(\varepsilon^2))^{\frac{1}{2}}.
\end{split}
\end{equation}

\noindent The result follows immediately from (\ref{eq:33}) and (\ref{eq:34}).
\end{proof}

It has been shown in \cite[Proposition 2.1]{AOST} that if $(e_n)$ is the unit vector basis of
$\ell_1$ and $(f_n)$ is a normalized subsymmetric basic sequence which is not equivalent to $(e_n)$
then $(e_n) >> (f_n)$.  Thus since the unit vector basis of $S$ is normalized and subsymmetric
we have that the unit vector basis of $\ell_1$ s.c. dominates the unit vector basis of $S$.  Thus
Proposition~\ref{prop:stillD} gives:

\begin{Rmk} \label{rmk:2star}
The unit vector basis of $\ell_2$ s.c. dominates the unit vector basis of $S_2.$
\end{Rmk}

\begin{Thm}
There exists an infinite dimensional subspace $Y$ of $D$ having a basis and $T \in \mathcal{L}(Y,D)$  such that
$T \not \in \mathbb{C} i_{Y \rightarrow D} + \mathcal{K}(Y,D)$.
\end{Thm}
\begin{proof}
We will refer to $(s_i)$ as the unit vector basis of the Schlumprecht space $S$ and $(\sqrt{s_i})$
as the unit vector basis of $S_2$.
Then apply
Theorem~\ref{thm:AOST_p} for $p=2$ and the spreading model for the unit vector basis of $D$.
By Proposition~\ref{prop:spreading} we have that $(\sqrt{s_i})$ is the spreading model of the
unit vector basis of
$D$.  By Remark~\ref{rmk:star} we have that 2 is in the Krivine set of $(\sqrt{s_i})$.
By Remark~\ref{rmk:2star} we have that the unit vector basis of $\ell_2$ s.c. dominates
$(\sqrt{s_i})$.
Thus by Theorem~\ref{thm:AOST_p} we have there exists an infinite dimensional subspace $Y$ of
$D$ and $T \in \mathcal{L}(Y,D)$ such that
$T \not \in \mathbb{C}i_{Y\rightarrow D} + \mathcal{K}(Y,D).$
\end{proof}

\end{document}